\documentclass{article}
\usepackage{amsmath}
\usepackage{graphicx}
\usepackage{amsfonts}
\usepackage{amssymb}
\newtheorem{theorem}{Theorem}
\newtheorem{acknowledgement}[theorem]{Acknowledgement}

\newtheorem{claim}[theorem]{Claim}

\newtheorem{corollary}[theorem]{Corollary}

\newtheorem{definition}[theorem]{Definition}

\newtheorem{lemma}[theorem]{Lemma}

\newtheorem{proposition}[theorem]{Proposition}
\newtheorem{remark}[theorem]{Remark}

\newenvironment{proof}[1][Proof]{\noindent\textbf{#1.} }{\ \rule{0.5em}{0.5em}}

\begin{document}

\title{Collapsing sequences of solutions to the Ricci flow on 3-manifolds with almost
nonnegative curvature}
\author{Bennett Chow\\University of California, San Diego
\and David Glickenstein\\University of California, San Diego
\and Peng Lu\\University of Oregon}
\date{}
\maketitle

\section{Introduction}

\label{result is general}\label{give generalization of Hamilton's classification}

We shall prove a general result about sequences of solutions to the Ricci flow
on compact or complete noncompact 3-manifolds with locally uniformly bounded
almost nonnegative sectional curvatures and diameters tending to infinity. It
is known that such sequences occur when one dilates about a singularity of a
solution to the Ricci flow on a 3-manifold. Our main result assumes collapse
and is complementary to the injectivity radius estimate in \cite[\S
25]{Hamilton Formation} and \cite{CKL}. In particular, the bump-like point
condition required in those papers is not assumed here; instead, we assume
collapse, which rules out bump-like points. As an application of our result we
give a generalization of Hamilton's singularity theory in dimension 3 to
classify collapsed singularity models arising from Type IIb (infinite time)
singularities. From Perelman's work \cite[\S\S4 and 7]{Pe}, such collapsed
singularity models cannot occur as limits of dilations about \emph{finite
time} singularities. Besides obtaining a local injectivity radius estimate and
ruling out the cigar as a limit for finite time singularities, Perelman has
also enlarged the class of points and times about which one can obtain good
limits. In view of Perelman's improvements of Hamilton's singularity theory,
one expects to be able to combine Perelman's ideas with Fukaya's ideas; this
is not discussed here.

\label{nice local limit geometries}

Limits of collapsing 3-manifolds with lower curvature bounds are Alexandrov
spaces with integer dimension 1 or 2 (we rule out dimension 0 by assuming the
diameters tend to infinity). In the study of singularities, the sequences of
solutions of the Ricci flow on 3-manifolds which arise have almost nonnegative
sectional curvatures. This assumption, together with the smoothing properties
of the Ricci flow (especially the strong maximum principle), put strong
restrictions on the local geometries of the solutions in the sequence. In
particular, the limit local covering geometries are locally the products of
positively curved surfaces with the real line.

\label{virtual 2d limit}\label{may satisfy modified RF}

It is partly for the above reason that we shall be able to extract a\emph{
virtual} 2-dimensional limit solution of the Ricci flow. The reason we call
this `limit' solution `virtual' is that it is not actually a limit of the
sequence, but rather constructed from limits of local covers of the sequence.
We expect that it is a limit of covers of exhaustions of the solutions in the sequence.

When the limit space is 2-dimensional, this allows us to extract a virtual
solution to the Ricci flow when the actual limit is an orbifold which may not
be a solution. The possible types of singularities of the orbifold are:
$D^{2}/\mathbb{Z}_{p}$ with rotation action, $D^{2}/\mathbb{Z}_{2}$ with
reflection action and $D^{2}/\mathbb{D}_{2p},$ where $D^{2}$ is the 2-disk and
$\mathbb{D}_{2p}$ is the dihedral group of order $2p$ for some $p>1$. We shall
show that in this case the orbifold limit has at most 1 singular point of type
$D^{2}/\mathbb{Z}_{p}$ or $D^{2}/\mathbb{D}_{2p}$. In particular, it is a good
orbifold with a finite cyclic or dihedral cover diffeomorphic to the plane.

The virtual limit associated to a 1-dimensional limit space is rotationally
symmetric, complete, noncompact, with bounded positive curvature. The
advantage of obtaining a 2-dimensional virtual limit as compared to a
1-dimensional actual limit is that 1-dimensional spaces have no intrinsic
geometry except for distances and in particular have no nontrivial curvature.

Our construction of the virtual limit relies on Hamilton's strong maximum
principle for systems (see \cite{Hamilton 4d PCO}), Fukaya's local covering
geometry theory (see \cite{F3}), and a Cheeger-Gromov type compactness theorem
for the Ricci flow (see \cite{Hamilton Compactness} and \cite{Gl}). The reader
is also directed to \cite{CM} for an application of Gromov-Hausdorff distance
in the study of the Ricci flow.\label{2d orbifold limit has only 1 singular point}

Two abbreviations we shall commonly use are GH for Gromov-Hausdorff and RF for
Ricci flow.

\begin{acknowledgement}
We would like to thank Peter Petersen, Nolan Wallach, and McKenzie Wang for
very helpful discussions.
\end{acknowledgement}

\section{The sequences of solutions with almost nonnegative sectional curvature}

\subsection{The definition of almost nonnegative sectional curvature}

Let $\left\{  \left(  M_{i}^{3},g_{i}\left(  t\right)  ,O_{i}\right)  :\text{
}t\in\left(  \alpha,\omega\right)  \right\}  _{i\in\mathbb{N}},$ where
$\alpha<0$ and $\omega>0,$ be a sequence of orientable complete solutions to
the RF with origins $O_{i}\in M_{i}^{3}$. Let $Rm_{i}\left(  x,t\right)  $
denote the Riemannian curvature operator of $g_{i}\left(  t\right)  .$

\begin{definition}
(i) The sequence $\left\{  \left(  M_{i}^{3},g_{i}\left(  t\right)
,O_{i}\right)  \right\}  $ is said to have \textbf{locally uniformly bounded
geometry} (or \textbf{bounded geometry} for short) if for every closed
subinterval $\left[  \beta,\psi\right]  \subset\left(  \alpha,\omega\right)
,$ $\rho<\infty,$ and $k\in\mathbb{N}\cup\left\{  0\right\}  $ there exists a
constant $C=C\left(  \beta,\psi,\rho,k\right)  <\infty$ such that
\[
\sup_{\ell\leq k,\,\,i\in\mathbb{N}}\left\{  \max_{B_{g_{i}\left(  0\right)
}\left(  O_{i},\rho\right)  \times\left[  \beta,\psi\right]  }\left|
\nabla^{\ell}Rm_{i}(x,t)\right|  _{g_{i}\left(  t\right)  }\right\}  \leq C.
\]

(ii) We say that the sequence $\left\{  \left(  M_{i}^{3},g_{i}\left(
t\right)  ,O_{i}\right)  \right\}  $ has \textbf{bounded diameters} if there
exists a constant $C<\infty$ such that
\[
\text{diam}\left(  M_{i}^{3},g_{i}\left(  0\right)  \right)  \leq C\text{, for
all }i\in\mathbb{N}.
\]
If
\[
\lim_{i\rightarrow\infty}\text{diam}\left(  M_{i}^{3},g_{i}\left(  0\right)
\right)  =\infty
\]
we say the sequence has \textbf{unbounded diameters}.

(iii) We say that the origins are \textbf{essential} if there exists $c>0$
such that $\left|  Rm_{i}\left(  O_{i},0\right)  \right|  _{g_{i}\left(
0\right)  }\geq c$ for all $i\in\mathbb{N}.$
\end{definition}

Let $\lambda_{1}\left(  Rm_{i}\right)  \leq\lambda_{2}\left(  Rm_{i}\right)
\leq\lambda_{3}\left(  Rm_{i}\right)  $ denote the eigenvalues of $Rm_{i}$.
Note that $\lambda_{k}\left(  Rm_{i}\right)  $ is twice the sectional curvature.

\begin{definition}
(i) The sequence is said to have \textbf{almost nonnegative sectional
curvatures} (or \textbf{ANSC}) if it has bounded geometry and for every
$\left[  \beta,\psi\right]  \subset\left(  \alpha,\omega\right)  $ and
$\rho<\infty,$ there exists $\delta_{i}\searrow0$ such that%
\[
\lambda_{1}\left(  Rm_{i}\right)  \geq-\delta_{i}\qquad\text{on }%
B_{g_{i}\left(  0\right)  }\left(  O_{i},\rho\right)  \times\left[  \beta
,\psi\right]  \text{ for all }i\in\mathbb{N}.
\]

(ii) An \textbf{injectivity radius estimate holds for the origins} if there
exists $c>0$ such that
\[
\left|  Rm_{i}\left(  O_{i},0\right)  \right|  _{g_{i}\left(  0\right)
}~\text{inj}_{g_{i}\left(  0\right)  }\left(  O_{i}\right)  ^{2}\geq
c\qquad\text{for all }i\in\mathbb{N}.
\]

(iii) The sequence is said to have \textbf{collapsing origins} if%
\[
\left|  Rm_{i}\left(  O_{i},0\right)  \right|  _{g_{i}\left(  0\right)
}~\text{inj}_{g_{i}\left(  0\right)  }\left(  O_{i}\right)  ^{2}%
\rightarrow0\qquad\text{as }i\rightarrow\infty.
\]
\end{definition}

By passing to a subsequence we may assume that one of the two alternatives
(ii) or (iii) holds.

\begin{definition}
An ANSC sequence $\left\{  \left(  M_{i}^{3},g_{i}\left(  t\right)
,O_{i}\right)  \right\}  $ has \textbf{bump-like origins} if there exists
$c>0$ such that
\[
\frac{\lambda_{1}\left(  Rm_{i}\right)  \left(  O_{i},0\right)  }{\left|
Rm_{i}\left(  O_{i},0\right)  \right|  _{g_{i}\left(  0\right)  }}\geq c
\]
for all $i\in\mathbb{N}.$ If $\lim_{i\rightarrow\infty}\lambda_{1}\left(
Rm_{i}\right)  \left(  O_{i},0\right)  /\left|  Rm_{i}\left(  O_{i},0\right)
\right|  _{g_{i}\left(  0\right)  }=0$ then we say the origins are
\textbf{split-like}.
\end{definition}

Throughout this paper we shall assume that the sequence $\left\{  \left(
M_{i}^{3},g_{i}\left(  t\right)  ,O_{i}\right)  \right\}  $ has bounded
geometry. In this section we further assume that the sequence has ANSC with
essential origins.

\subsection{The four types of sequences}

We shall categorize the types of sequences by whether they have bump-like or
split-like origins and whether the diameters are bounded or unbounded.

\textbf{T1. The essential bump-like origins with bounded diameters case}. In
this case we assume that the sequence $\left\{  \left(  M_{i}^{3},g_{i}\left(
t\right)  ,O_{i}\right)  \right\}  $ has essential bump-like origins and
bounded diameter. Then a subsequence of $M_{i}^{3}$ are diffeomorphic to
spherical space forms.

\begin{lemma}
Let $\left\{  \left(  M_{i}^{3},g_{i}\left(  t\right)  ,O_{i}\right)
\right\}  $ be a sequence of solutions to the RF with ANSC, essential
bump-like origins, and bounded diameters. Then there exists a subsequence such
that $g_{i}\left(  0\right)  $ has positive sectional curvature and $M_{i}%
^{3}$ is diffeomorphic to a spherical space form for all $i.$
\end{lemma}

\begin{proof}
It follows from the proof of Lemma 25.2 in \cite{Hamilton Formation}, using
the bounded diameters assumption, that we get positive sectional curvature on
all of $M_{i}^{3}$ for a subsequence. The lemma now follows from Hamilton's
classification of compact 3-manifolds with positive Ricci curvature
\cite{Hamilton 3-manifold}.
\end{proof}

\textbf{T2. The essential bump-like origins with unbounded diameters case -
Hamilton's injectivity radius estimate}. In this case we assume that the
sequence $\left\{  \left(  M_{i}^{3},g_{i}\left(  t\right)  ,O_{i}\right)
\right\}  $ has essential bump-like origins and unbounded diameter. Then there
is an injectivity radius estimate for a subsequence and hence there is a
subsequence which converges to a complete solution diffeomorphic to
$\mathbb{R}^{3}$ with positive sectional curvature.

\begin{proposition}
Let $\left\{  \left(  M_{i}^{3},g_{i}\left(  t\right)  ,O_{i}\right)
\right\}  $ be a sequence of solutions to the RF with ANSC, essential
bump-like origins, and unbounded diameters. Then there exists a subsequence
and a constant $c>0$ such that%
\[
\text{inj}_{g_{i}\left(  0\right)  }\left(  O_{i}\right)  \geq c.
\]
\end{proposition}

For the proof see \cite[\S25]{Hamilton Formation} and \cite{CKL}.

\textbf{T3. The essential split-like origins with bounded diameters case}. In
this case we assume that the sequence $\left\{  \left(  M_{i}^{3},g_{i}\left(
t\right)  ,O_{i}\right)  \right\}  $ has essential split-like origins and
bounded diameters. If the origins are not collapsing, then there is an
injectivity radius estimate for a subsequence and hence there exists a further
subsequence which converges to a\ compact $\left\{  \left(  M_{\infty}%
^{3},g_{\infty}\left(  t\right)  ,O_{\infty}\right)  \right\}  $ with
$\lambda_{1}\left(  Rm_{g_{\infty}\left(  t\right)  }\right)  =\lambda
_{2}\left(  Rm_{g_{\infty}\left(  t\right)  }\right)  =0$ and $\lambda
_{3}\left(  Rm_{g_{\infty}\left(  t\right)  }\right)  >0.$ Such solutions are
classified in \cite{Hamilton 4d PCO} using the strong maximal principle;
$M_{\infty}^{3}$ is diffeomorphic to $S^{2}\times S^{1}$ or the twisted
product $S^{2}\tilde{\times}S^{1}$. If the origins are collapsing there is a
subsequence such that $M_{i}^{3}$ are graph manifolds.

\begin{lemma}
Let $\left\{  \left(  M_{i}^{3},g_{i}\left(  t\right)  ,O_{i}\right)
\right\}  $ be a sequence of solutions to the RF with ANSC, essential
collapsing origins, and bounded diameters. Then for any $\varepsilon>0$ there
exists a subsequence such that $M_{i}^{3}$ is $\varepsilon$-collapsed for all
$i,$ i.e. the supremum of the injectivity radii of $M_{i}^{3}$ is at most
$\varepsilon$. When $\varepsilon$ is small enough each $M_{i}^{3}$ is
diffeomorphic to a graph manifold.
\end{lemma}

\begin{proof}
By the injectivity radius decay estimate of \cite{CGT} or \cite{CLY} and the
bounded diameter assumption, the sequence of Riemannian manifolds $\left\{
\left(  M_{i}^{3},g_{i}\left(  0\right)  \right)  \right\}  $ collapses. That
is, $\max_{x\in M_{i}}$inj$_{g_{i}(0)}\left(  x\right)  \rightarrow0$ as
$i\rightarrow\infty.$ The topological conclusion follows from Cheeger-Gromov
theory (see \cite{CG1}, \cite{CG2}), as there exists an F-structure on $M_{i}$
for $i$ large enough. For more details see the discussion in \cite[p.548]{R}.
\end{proof}

\textbf{T4. The essential split-like origins with unbounded diameters case}.
In this case we assume that the sequence $\left\{  \left(  M_{i}^{3}%
,g_{i}\left(  t\right)  ,O_{i}\right)  \right\}  $ has essential split-like
origins and unbounded diameters. If the origins are not collapsing, then there
is an injectivity radius estimate at $O_{i}$ for a subsequence and there
exists a subsequence which converges to a\ complete noncompact $\left\{
\left(  M_{\infty}^{3},g_{\infty}\left(  t\right)  ,O_{\infty}\right)
\right\}  $ with $\lambda_{1}\left(  Rm_{g_{\infty}\left(  t\right)  }\right)
=\lambda_{2}\left(  Rm_{g_{\infty}\left(  t\right)  }\right)  =0$ and
$\lambda_{3}\left(  Rm_{g_{\infty}\left(  t\right)  }\right)  >0.$ Such
solutions can be classified using the strong maximal principle as in
\cite{Hamilton 4d PCO}; the universal cover of $M_{\infty}^{3}$ is
diffeomorphic to $\Sigma^{2}\times\mathbb{R}^{1}$ where $\Sigma^{2}$ is a
complete surface with positive curvature.

The main focus of this paper is on the case when $\left\{  \left(  M_{i}%
^{3},g_{i}\left(  t\right)  ,O_{i}\right)  \right\}  $ has ANSC, unbounded
diameters, and essential collapsing split-like origins.

\subsection{Dilations about sequences of points and times\label{sec 2.3
dilate}}

In this subsection we show how ANSC sequences $\left\{  \left(  M_{i}%
^{3},g_{i}\left(  t\right)  ,O_{i}\right)  \right\}  $ can be possibly
obtained by dilation of the RF at singularities. Recall the best estimate of
the sectional curvatures tending to nonnegative due to Hamilton \cite[Theorem
4.1]{Hamilton Nonsingular}, which improves earlier estimates of \cite[Theorem
24.4]{Hamilton Formation} and \cite[Theorem 2]{Ivey}.

\begin{proposition}
\label{Prop Curvature Tends to Nonnegative}Let $\left(  M^{3},g\left(
t\right)  \right)  ,$ $t\in\lbrack0,T),$ be a complete solution to the RF with
time-dependent bounded curvature, i.e. sup$_{x\in M^{3}}\left|  Rm\left(
x,t\right)  \right|  _{g\left(  t\right)  }\leq C\left(  t\right)  $. If
$\lambda_{1}\doteqdot\lambda_{1}\left(  Rm\right)  \geq-C_{0}$ at time $0$ for
some $C_{0}>0,$ then at any point and time where $\lambda_{1}<0$ we have%
\[
R\geq-\lambda_{1}\left[  \ln\left(  -\lambda_{1}\right)  +\ln\left(
1+C_{0}t\right)  -\ln C_{0}-3\right]  ,
\]
where $R$ is the scalar curvature.
\end{proposition}

Let $\left(  M^{3},g\left(  t\right)  \right)  ,$ $t\in\lbrack0,\infty),$ be a
solution to the RF. Let $\left(  x_{i},t_{i}\right)  $ be a sequence of
spacetime points and $K_{i}\doteqdot\left\vert Rm\left(  x_{i},t_{i}\right)
\right\vert _{g\left(  t_{i}\right)  }.$ We say that $\left(  x_{i}%
,t_{i}\right)  $ is \textbf{dilatable} if there exists $\beta<0$ and $\psi>0$
such that for every $\rho<\infty$ there exists $C<\infty$ such that
\[
\left\vert Rm\left(  x,t\right)  \right\vert _{g\left(  t\right)  }\leq CK_{i}%
\]
for $x\in B_{g\left(  t_{i}\right)  }\left(  x_{i},\rho/\sqrt{K_{i}}\right)  $
and $t\in\left[  t_{i}+\beta/K_{i},t_{i}+\psi/K_{i}\right]  .$ Let
$g_{i}\left(  t\right)  \doteqdot K_{i}g(t_{i}+t/K_{i})$. If $\left(
x_{i},t_{i}\right)  $ is dilatable, then $\left\vert Rm_{i}\left(  x,t\right)
\right\vert _{g_{i}\left(  t\right)  }\leq C$ on $B_{g_{i}\left(  0\right)
}\left(  x_{i},\rho\right)  \times\left[  \beta,\psi\right]  .$

A sequence $\left(  x_{i},t_{i}\right)  $ with $t_{i}\rightarrow\infty$ is
called \textbf{Type III-like} if there exists a constant $C<\infty$ such that
$K_{i}\leq C/t_{i}$ for all $i\in\mathbb{N}.$ If $\lim_{i\rightarrow\infty
}t_{i}K_{i}=\infty$ the sequence $\left(  x_{i},t_{i}\right)  $ is said to be
\textbf{Type IIb-like}. Let $\lambda_{k}\left(  i\right)  \doteqdot\lambda
_{k}\left(  Rm_{i}\right)  \left(  x_{i},t_{i}\right)  $ for $k=1,2,3$. The
sequence $\left(  x_{i},t_{i}\right)  $ has \textbf{almost nonnegative
sectional curvatures} if $\lambda_{3}\left(  i\right)  >0$ for $i$ large
enough and $\lim_{i\rightarrow\infty}\min\left\{  \lambda_{1}\left(  i\right)
,0\right\}  /\lambda_{3}\left(  i\right)  =0.$

\begin{corollary}
Let $\left\{  \left(  x_{i},t_{i}\right)  \right\}  $ be a dilatable and Type
IIb-like sequence, then $\left(  x_{i},t_{i}\right)  $ has almost nonnegative
sectional curvatures and the dilated solutions $\{(M^{3},g_{i}\left(
t\right)  ,x_{i})\},$ $t\in\left[  \beta,\psi\right]  $ have ANSC.
\end{corollary}

\begin{proof}
If $\lambda_{1}\left(  x,t_{i}\right)  <0$ where $x\in B_{g\left(
t_{i}\right)  }\left(  x_{i},\rho/\sqrt{K_{i}}\right)  ,$ then by Proposition
\ref{Prop Curvature Tends to Nonnegative} and $R\left(  x,t_{i}\right)
=\lambda_{1}\left(  x,t_{i}\right)  +\lambda_{2}\left(  x,t_{i}\right)
+\lambda_{3}\left(  x,t_{i}\right)  ,$%
\begin{equation}
\lambda_{3}\left(  x,t_{i}\right)  \geq-\frac{1}{2}\lambda_{1}\left(
x,t_{i}\right)  \ln\left[  -\lambda_{1}\left(  x,t_{i}\right)  \left(
C_{0}^{-1}+t_{i}\right)  e^{-2}\right]  . \label{scalar curv inequality}%
\end{equation}
First consider $x=x_{i}.$ For any $L\in(1,\infty),$ if $-\lambda_{1}\left(
i\right)  >e^{2L+2}\left(  C_{0}^{-1}+t_{i}\right)  ^{-1},$ then $\lambda
_{3}\left(  i\right)  /\left[  -\lambda_{1}\left(  i\right)  \right]  \geq
L>1$ by (\ref{scalar curv inequality}). Note that this implies that
$\lambda_{3}^{2}\left(  i\right)  >\lambda_{1}^{2}\left(  i\right)  .$ If
$0\leq-\lambda_{1}\left(  i\right)  \leq e^{2L+2}\left(  C_{0}^{-1}%
+t_{i}\right)  ^{-1}$ then, by the Type IIb assumption, for $i$ large enough
\begin{equation}
K_{i}=\left(  \sum_{k=1}^{3}\lambda_{k}\left(  i\right)  ^{2}\right)
^{1/2}\geq\sqrt{3}Le^{2L+2}\left(  C_{0}^{-1}+t_{i}\right)  ^{-1}\geq-\sqrt
{3}L\lambda_{1}\left(  i\right)  . \label{K is like lambda1}%
\end{equation}
This implies (since $L>1$) that $\lambda_{2}^{2}\left(  i\right)  +\lambda
_{3}^{2}\left(  i\right)  >2\lambda_{1}^{2}\left(  i\right)  .$ Therefore
$\lambda_{3}^{2}\left(  i\right)  >\lambda_{1}^{2}\left(  i\right)  $ with no
condition on $\lambda_{1}\left(  i\right)  ,$ so $\sqrt{3}\lambda_{3}\left(
i\right)  \geq K_{i}\geq\lambda_{3}\left(  i\right)  .$ Hence, by (\ref{K is
like lambda1}), $\sqrt{3}\lambda_{3}\left(  i\right)  \geq K_{i}\geq-\sqrt
{3}L\lambda_{1}\left(  i\right)  ,$ or $\lambda_{3}\left(  i\right)  /\left[
-\lambda_{1}\left(  i\right)  \right]  \geq L.$ Since $L$ is arbitrary, we
have that $\left(  x_{i},t_{i}\right)  $ has almost nonnegative sectional curvatures.

Now consider any $x\in B_{g\left(  t_{i}\right)  }\left(  x_{i},\rho
/\sqrt{K_{i}}\right)  .$ We want to show that $g_{i}\left(  t\right)  $ has
ANSC, so we need to show that $\lambda_{1}\left(  g_{i}\left(  t\right)
\right)  =\lambda_{1}\left(  x,t\right)  /K_{i}\geq-\delta_{i}.$ Recall that
$\lambda_{3}\left(  i\right)  $ is comparable to $K_{i},$ so it is sufficient
to prove $\lambda_{1}\left(  x,t\right)  /\lambda_{3}\left(  i\right)
\geq-\delta_{i}$. For any $L\in(1,\infty),$ if $-\lambda_{1}\left(
x,t\right)  >e^{2L+2}\left(  C_{0}^{-1}+t\right)  ^{-1},$ then $\lambda
_{3}\left(  x,t\right)  /\left[  -\lambda_{1}\left(  x,t\right)  \right]  \geq
L$ by (\ref{scalar curv inequality}). This implies $\lambda_{3}\left(
x,t\right)  $ is comparable to $K\left(  x,t\right)  ,$ i.e.
\[
\lambda_{3}\left(  x,t\right)  \leq K\left(  x,t\right)  \leq\sqrt{3}%
\lambda_{3}\left(  x,t\right)  ,
\]
as above. By the dilatable assumption $C\left(  \rho\right)  K_{i}\geq
K\left(  x,t\right)  ,$ so $C\left(  \rho\right)  \sqrt{3}\lambda_{3}\left(
i\right)  \geq\lambda_{3}\left(  x,t\right)  $ and $\lambda_{3}\left(
i\right)  /\left[  -\lambda_{1}\left(  x,t\right)  \right]  \geq\left(
C\left(  \rho\right)  \sqrt{3}\right)  ^{-1}L.$

Now suppose $0\leq-\lambda_{1}\left(  x,t\right)  \leq e^{2L+2}\left(
C_{0}^{-1}+t\right)  ^{-1}.$ Since $t\in\left[  t_{i}+\beta/K_{i},t_{i}%
+\psi/K_{i}\right]  $ and $\lim t_{i}K_{i}=\infty,$ we get for all $t$ that
$t/t_{i}\geq1+\frac{\beta}{t_{i}K_{i}}\geq1/2$ for $i$ large enough. Because
$\sqrt{3}\lambda_{3}\left(  i\right)  \geq K_{i},$ we have, for $i$ large
enough (using the Type IIb assumption),
\[
\lambda_{3}\left(  i\right)  \geq\frac{1}{\sqrt{3}}K_{i}\geq Le^{2L+2}\left(
C_{0}^{-1}+\frac{1}{2}t_{i}\right)  ^{-1}\geq Le^{2L+2}\left(  C_{0}%
^{-1}+t\right)  ^{-1}\geq-L\lambda_{1}\left(  x,t\right)  .
\]
We have proven that for any $L\in(1,\infty)$ and for $i$ large enough we have
$\lambda_{3}\left(  i\right)  /\left[  -\lambda_{1}\left(  x,t\right)
\right]  \geq C^{\prime}\left(  \rho\right)  ^{-1}L.$ We conclude that
$g_{i}\left(  t\right)  $ has ANSC.
\end{proof}

The about result is not true for Type III-like sequences of points; for
example, constant negative sectional curvature solutions are of Type III.

\section{Review of Fukaya's local theory}

One of the main tools we shall use is Fukaya's local theory, which describes
the local geometry of collapsed limits. In dimension 3 the types of the local
geometries are quite limited. In this section when we use the ANSC assumption
we shall make it explicit.

\subsection{Fukaya's main theorem\label{B1b.}}

Recall Definition 0-4 in \cite{F3}.\footnote{We use the terminology
\emph{nice} instead of Fukaya's \emph{smooth} to distinguish it from
$C^{\infty}.$}\label{definition of nice (Fukaya's smooth)}

\begin{definition}
We say that a metric space $\left(  X,d\right)  $ is \textbf{nice} if for
every point $p\in X,$ there exists

\begin{enumerate}
\item  a neighborhood $U$ of $p$ in $X,$ and a neighborhood $V$ of $\vec{0}$
in $\mathbb{R}^{m}$ for some $m\in\mathbb{N}\cup\left\{  0\right\}  ,$

\item  a compact Lie group $\Gamma$ with a faithful representation of $\Gamma$
into $O\left(  m,\mathbb{R}\right)  $ where $\Gamma^{0}$ (the identity
component of $\Gamma$) is isomorphic to a torus, and

\item  a $\Gamma$-invariant Riemannian metric $h$ on $V,$
\end{enumerate}
\end{definition}

\noindent\emph{such that }$\left(  U,\left.  d\right|  _{U}\right)  $\emph{ is
isometric to }$\left.  \left(  V,d_{h}\right)  \right/  \Gamma,$\emph{ which
is a metric space with distance function induced on the quotient }$V/\Gamma
$\emph{ from }$h$\emph{ (so that, in particular, }$U$\emph{ is homeomorphic to
}$V/\Gamma$\emph{).}

\label{bounded geom implies limit is nice}In \cite{F3} Fukaya proves that
given a sequence of Riemannian manifolds with $C^{\infty}$-uniformly bounded
geometry there is a subsequence which converges to a nice metric space. This
implies that nice metric spaces are Lipschitz dense in the closure of
Riemannian manifolds of dimension $n$ with bounded curvature. We will outline
Fukaya's proof since we shall need elements of it in our classification of the
local geometries. In addition, we shall clarify the exact versions of Fukaya's
results which we will need and generalize them to the case of solutions of the
RF as done in \cite{Gl}.

\subsection{Construction of the limit metric, local group $G_{\infty}%
$\label{sec 3.2 limit metric}}

Recall the definition of local groups from \cite[23D]{P}; these are sometimes
called pseudogroups.\label{local group def} For example, a neighborhood of the
identity in a Lie group is a local group. It will be important to consider
local groups which are not connected.

Let $\left\{  \left(  M_{i}^{3},g_{i}\left(  t\right)  ,O_{i}\right)
:t\in\left(  \alpha,\omega\right)  \right\}  $ be a sequence of solutions of
the RF. Suppose $\left(  M_{i}^{3},g_{i}\left(  0\right)  ,O_{i}\right)  $
converges to the metric space $\left(  X_{\infty},d_{\infty}(0),O_{\infty
}\right)  $ in the pointed GH topology. Fix an $\epsilon>0$ and $P_{\infty}\in
X_{\infty}$ and let $P_{i}\in M_{i}^{3}$ such that $B_{g_{i}(0)}^{3}\left(
P_{i},1+\epsilon\right)  $ converges to $B_{d_{\infty}(0)}^{3}\left(
P_{\infty},1+\epsilon\right)  \subset X_{\infty}$ in GH. Fix a frame $F_{i}$
of $T_{P_{i}}M_{i}^{3}$ orthonormal with respect to the metric $g_{i}(0);$
these frames allows us to identify each unit ball in $T_{P_{i}}M_{i}^{3}$
centered at the origin in $T_{P_{i}}M_{i}^{3}$ with the Euclidean unit ball
$B^{3}\left(  1\right)  \subset\mathbb{R}^{3}$ centered at the origin $\vec
{0}.$ Consider the exponential map for $g_{i}\left(  0\right)  $ restricted to
ball:
\begin{equation}
\exp_{P_{i}}\doteqdot\left.  \exp_{P_{i}}^{g_{i}(0)}\right|  _{B^{3}\left(
1+\epsilon\right)  }:B^{3}(1+\epsilon)\rightarrow M_{i}^{3}.\nonumber
\end{equation}
Assume that the sectional curvatures $K_{g_{i}(0)}\leq1$ on the ball
$B_{g_{i}(0)}\left(  P_{i},1+\epsilon\right)  $ so that $\exp_{P_{i}}$ is a
local diffeomorphism. We consider the pulled-back metrics $\widetilde{g}%
_{i}(t)\doteqdot\exp_{P_{i}}^{\ast}g_{i}(t)$ on $B^{3}\left(  1+\epsilon
\right)  $ for $t\in\left(  \alpha,\omega\right)  .$ By the Arzela-Ascoli
theorem and the bounded geometry assumption, there exists a subsequence (we
still denote the subsequence by $\widetilde{g}_{i}(t)$ and continue to use
this convention with further subsequences below) such that $\widetilde{g}%
_{i}(t)$ converges in $C^{\infty}$ on $B^{3}\left(  1\right)  $ to a smooth
solution $\widetilde{g}_{\infty}(t)$ to the RF on $B^{3}(1).$ By a
diagonalization argument we may assume that this convergence holds for all
$P_{\infty}$ in a countable dense subset of $X_{\infty}.$

\label{limit isom group}We shall consider the set of continuous maps
$C^{0}(B^{3}\left(  1/2\right)  ,B^{3}(1))$ as a metric space with the metric
\[
d_{C}\left(  \gamma,\gamma^{\prime}\right)  =\sup_{x\in B^{3}\left(
1/2\right)  }d_{\tilde{g}_{\infty}\left(  0\right)  }\left[  \gamma\left(
x\right)  ,\gamma^{\prime}\left(  x\right)  \right]  .
\]
We define the sets $G_{i}$ which consist of local deck transformations of the
local covering map $\exp_{P_{i}}:B^{3}\left(  1\right)  \rightarrow M_{i}$ as
\[
G_{i}\doteqdot\left\{  \gamma\in C^{0}(B^{3}\left(  1/2\right)  ,B^{3}%
(1)):\exp_{P_{i}}\circ\,\gamma=\exp_{P_{i}}\right\}  .
\]
Clearly each $G_{i}$ is a discrete local group of local isometries of the
Riemannian manifold $\left(  B^{3}\left(  1\right)  ,\tilde{g}_{i}\left(
t\right)  \right)  $ for each $t$. Gromov \cite{Gromov} calls this group the
local fundamental pseudogroup and there is a geometric description of elements
of $G_{i}$ in \cite[\S4]{CGT}.

\label{G_infty as Haus limit of G_i}We will define a limit group $G_{\infty}$.
Since $G_{i}$ are local isometries of $\left(  B^{3}\left(  1\right)
,\tilde{g}_{i}\left(  0\right)  \right)  $ and $\tilde{g}_{i}\left(  0\right)
$ converge to $\tilde{g}_{\infty}\left(  0\right)  $ on $B^{3}\left(
1\right)  $ in $C^{\infty},$ for large $i$ each $G_{i}$ is a closed subset of
the following set of \emph{quasi-isometries}\footnote{A quasi-isometry is a
homeomorphism which distorts distances by a bounded amount.}%
\[
L\doteqdot\left\{  \gamma\in C^{0}(B^{3}\left(  1/2\right)  ,B^{3}%
(1)):\frac{1}{2}\leq\frac{d_{\tilde{g}_{\infty}\left(  0\right)  }\left[
\gamma\left(  x\right)  ,\gamma\left(  y\right)  \right]  }{d_{\tilde
{g}_{\infty}\left(  0\right)  }\left[  x,y\right]  }\leq2\right\}  .
\]
By the Arzela-Ascoli theorem $L$ is compact. We consider the space
$\mathcal{S}$ of closed subsets of $L$ with the Hausdorff topology.
$\mathcal{S}$ is compact because $L$ is compact (see \cite[Theorem 7.3.8]%
{BBI}); hence there is a subsequence of $\left\{  G_{i}\right\}
\subset\mathcal{S}$ which converges to a set $G_{\infty}\in\mathcal{S}.$

The inclusion $G_{\infty}\subset C^{0}\left(  B^{3}\left(  1/2\right)
,B^{3}\left(  1\right)  \right)  $ defines a local action of $G_{\infty}$ on
$B^{3}\left(  1/2\right)  $. Elements of $G_{\infty}$ are local isometries
with respect to $\tilde{g}_{\infty}\left(  0\right)  ,$ that is, if $\gamma\in
G_{\infty}$ then $\gamma^{\ast}\left(  \tilde{g}_{\infty}\left(  0\right)
\right)  =\left.  \tilde{g}_{\infty}\left(  0\right)  \right|  _{B^{3}\left(
1/2\right)  }.$ The product operation in $G_{i}$ gives a product operation in
the limit $G_{\infty}$ which makes it a local group.

$G_{\infty}$ is a Lie group germ acting smoothly \cite[Lemma 3.1]{F3}. A Lie
group germ is a local group isomorphic to a neighborhood of the identity of a
Lie group. Furthermore, $G_{\infty}$ is nilpotent and has a neighborhood of
the identity where the exponential map of the Lie algebra is onto \cite[Lemma
4.1]{F3}.

\label{identity component of G_infty is abelian}Let $H^{0}$ denote the
identity component of a local group $H$. We have the following important
consequence of the nilpotency of Lie($G_{\infty}$) together with $\dim
M_{i}=3$. We shall use some results from \S\ref{Section 4 limit space}.

\begin{lemma}
\label{G abelian} Suppose the GH limit $X_{\infty}$ of an essential ANSC
sequence has Hausdorff $\dim X_{\infty}\geq1$, then $G_{\infty}^{0}$ is abelian.
\end{lemma}

\begin{proof}
It follows from the ANSC assumption and \S\ref{sec 4.2 conse strong max} that
$G_{\infty}$ acts isometrically on $D^{2}\left(  1/2\right)  \times(-1/2,1/2)$
with metric $^{2}\widetilde{h}_{\infty}(t)+du^{2}$. Let $Y$ be a vector field
on $D^{2}\left(  1/2\right)  \times(-1/2,1/2)$ generated by the $G_{\infty
}^{0}$ action. Then by Lemma\ \ref{structure of group action} we have
$Y=K+b\frac{\partial}{\partial u}$ where $K$ is a Killing vector field on
$D^{2}\left(  1/2\right)  $. Let $Y_{1}$ and $Y_{2}$ be two such vector fields
on $D^{2}\left(  1/2\right)  \times(-1/2,1/2)$, then $[Y_{1},Y_{2}%
]=[K_{1},K_{2}]$ and hence we get a Lie algebra Lie$_{0}=\left\{
K:K+b\frac{\partial}{\partial u}\in\text{Lie}\left(  G_{\infty}^{0}\right)
\right\}  $ of Killing vector fields on $D^{2}\left(  1/2\right)  $. By Lemma
4.1 of \cite{F3}, Lie$\left(  G_{\infty}\right)  =$Lie$\left(  G_{\infty}%
^{0}\right)  $ is nilpotent, this implies that Lie algebra Lie$_{0}$ is
nilpotent. Since Lie$_{0}$ consists of Killing vector fields on $D^{2}\left(
1/2\right)  $ with positive curvature, Lie$_{0}$ must be one dimensional.
Hence $G_{\infty}^{0}$ is an abelian local Lie group.
\end{proof}

\begin{remark}
We think this can be proven without the ANSC assumption using the results of
\cite[\S5]{F3}.
\end{remark}

Next we will prove that $\dim\left[  G_{\infty}\left(  \vec{0}\right)
\right]  \geq1$ when $\dim X_{\infty}\leq2$, where $G_{\infty}\left(  \vec
{0}\right)  =\left\{  \gamma\left(  \vec{0}\right)  :\gamma\in G_{\infty
}\right\}  .$ The fact that $G_{\infty}\left(  \vec{0}\right)  $ is infinite
is stated in \cite{SY}. Let $\gamma_{i}$ be the minimizer of $d_{\tilde{g}%
_{i}\left(  0\right)  }\left[  \gamma\left(  \vec{0}\right)  ,\vec{0}\right]
$ over all $\gamma\in G_{i}.$ This corresponds to the shortest geodesic loop
in $\left(  M_{i}^{3},g_{i}\left(  0\right)  \right)  $ based at $P_{i}.$

\begin{claim}
\label{outside ball}Let $\gamma_{i}\in G_{i}$ as defined above, then for all
$\varepsilon,$ $0<\varepsilon<1/4,$ there exists $n_{i}$ such that $\gamma
_{i}^{n_{i}}$ exists and $d_{\tilde{g}_{i}\left(  0\right)  }\left[
\gamma_{i}^{n_{i}}\left(  \vec{0}\right)  ,\vec{0}\right]  >\varepsilon$ and
$d_{\tilde{g}_{i}\left(  0\right)  }\left[  \gamma_{i}^{n_{i}-1}\left(
\vec{0}\right)  ,\vec{0}\right]  \leq\varepsilon$
\end{claim}

\begin{proof}
Suppose there was $\varepsilon>0$ such that this is not true. Then
$d_{\tilde{g}_{i}\left(  0\right)  }\left[  \gamma_{i}^{n}\left(  \vec
{0}\right)  ,\vec{0}\right]  <\varepsilon$ for all $n>0.$ Hence $\gamma
_{i}^{n}$ exists for all $n.$ It follows from \cite[Lemma 4.6]{CGT} that for
any $k$, $\gamma_{i}^{k}\neq id$. Since $G_{i}$ is a closed subset of the
compact set $L,$ there is a subsequence $n_{j}$ which converges to $\gamma
_{i}^{\infty}\in$ $G_{i}$. But then we can find $n_{j}$ such that $\gamma
_{i}^{n_{j}}$ is arbitrarily close to $\gamma_{i}^{\infty}.$ This means that
$d_{\tilde{g}_{i}\left(  0\right)  }\left[  \left(  \gamma_{i}^{\infty
}\right)  ^{-1}\gamma_{i}^{n_{j}}\left(  \vec{0}\right)  ,\vec{0}\right]
\rightarrow0$ as $j\rightarrow\infty.$ Hence $\left(  M_{i}^{3},g_{i}%
(0)\right)  $ has arbitrarily small geodesic 1-gons based at $P_{i},$ a contradiction.
\end{proof}

\begin{lemma}
\label{group infinite}Suppose the GH limit $X_{\infty}$ has (Hausdorff) $\dim
X_{\infty}\leq2$. Then $\dim\left[  G_{\infty}\left(  \vec{0}\right)  \right]
\geq1.$
\end{lemma}

\begin{proof}
We shall construct a one-dimensional path in $G_{\infty}\left(  \vec
{0}\right)  .$ Since $\dim X_{\infty}\leq2$ and $d_{\tilde{g}_{i}\left(
0\right)  }\left[  \gamma_{i}\left(  \vec{0}\right)  ,\vec{0}\right]
=2\left[  inj_{P_{i}}\left(  g_{i}\left(  0\right)  \right)  \right]  $, we
must have $d_{\tilde{g}_{i}\left(  0\right)  }\left[  \gamma_{i}\left(
\vec{0}\right)  ,\vec{0}\right]  \rightarrow0$ as $i\rightarrow\infty.$ By
Claim \ref{outside ball}, there exists a subsequence $\left\{  \gamma_{i_{j}%
}\right\}  $ such that $\gamma_{i_{j}}\rightarrow\gamma_{\infty}$ with
$\gamma_{\infty}\left(  \vec{0}\right)  =\vec{0}$ and there exist $n_{i_{j}}$
such that $\gamma_{i_{j}}^{n_{i_{j}}}\rightarrow\gamma_{\infty}^{\prime}\in
G_{\infty}$ as $j\rightarrow\infty$ with $\gamma_{\infty}^{\prime}\left(
\vec{0}\right)  \neq\vec{0}.$ We can now take the path $\gamma_{t}%
\doteqdot\lim_{j\rightarrow\infty}\gamma_{i_{j}}^{\left\lfloor tn_{i_{j}%
}\right\rfloor },$ where $\left\lfloor q\right\rfloor $ is the greatest
integer less than or equal to $q,$ for $t\in\left[  0,1\right]  .$ This is a
path from $\gamma_{\infty}$ to $\gamma_{\infty}^{\prime},$ and hence
$\gamma_{t}\left(  \vec{0}\right)  $ is a path from $\vec{0}$ to
$\gamma_{\infty}^{\prime}\left(  \vec{0}\right)  $ inside the orbit
$G_{\infty}\left(  \vec{0}\right)  .$
\end{proof}

We define
\[
\Gamma\doteqdot\left\{  \gamma\in G_{\infty}:\gamma\left(  \vec{0}\right)
=\vec{0}\right\}  ,
\]
the isotropy sub-local group of $\vec{0}$ in $G_{\infty}.$ We have a faithful
representation of $\Gamma$ into $O\left(  3\right)  $ with metric $\tilde
{g}_{\infty}\left(  0\right)  $ at $\vec{0}$ defined by $\gamma\mapsto
\gamma_{\ast}\left(  \vec{0}\right)  :T_{\vec{0}}B^{3}\left(  1/4\right)
\rightarrow T_{\vec{0}}B^{3}\left(  1/4\right)  .$ Actually, the elements of
$\Gamma$ are orientation preserving because $M_{i}$ are orientable, so the
representation is into $SO\left(  3\right)  $. The representation is faithful
since an isometry is uniquely determined by its derivative at one point.

The action can be thought of as the linear action of $SO\left(  3\right)  .$
That is, every element of $\Gamma$ is the restriction to $B^{3}\left(
1/2\right)  $ of an element of $SO\left(  3\right)  .$ $\Gamma$ acts by
isometries on the Riemannian manifold $\left(  B^{3}\left(  1\right)
,g_{\mathbb{E}}\right)  $, where $g_{\mathbb{E}}$ is the Euclidean metric.
This follows by differentiating the action of $\Gamma$ since rays from the
origin are geodesics of both the pulled back metric and the Euclidean metric.

\subsection{Decomposition of the actions\label{sec 3.3 local model}}

Given a set $G$ acting on a metric space $\left(  X,d\right)  $, for $p\in X$
and $\varepsilon>0$ we define
\[
G\left(  p,\varepsilon\right)  \doteqdot\left\{  \gamma\in G:d\left(
\gamma\left(  p\right)  ,p\right)  <\varepsilon\right\}  .
\]
For any $\varepsilon_{1}\in(0,1/4],$ since $G_{i}$ act on $B^{3}\left(
1/2\right)  $ as isometries, there is a well defined equivalence relation
$\sim$ on $B^{3}\left(  \varepsilon_{1}\right)  $ defined by $x\sim y$ if and
only if there exists a $\gamma\in G_{i}\left(  \vec{0},2\varepsilon
_{1}\right)  $ such that $\gamma\left(  x\right)  =y$. Hence we have a
quotient $\left.  B^{3}\left(  \varepsilon_{1}\right)  \right/  G_{i}\left(
\vec{0},2\varepsilon_{1}\right)  $. Similarly we can define the quotient
$B^{3}\left(  \varepsilon_{1}\right)  /G_{\infty}(\vec{0},2\varepsilon_{1}).$
The metric space $B_{d_{\infty}\left(  0\right)  }^{X_{\infty}}\left(
P_{\infty},\varepsilon_{1}\right)  $ is isometric to $B^{3}\left(
\varepsilon_{1}\right)  /G_{\infty}(\vec{0},2\varepsilon_{1})$ with the
quotient distance induced by $\tilde{g}_{\infty}\left(  0\right)  $, which can
be seen using an equivariant version of pointed GH convergence (see \cite[p.
10]{F3}).

\label{The local geometry of the Alexandrov space limit}We would like to write
the group as $G_{\infty}\left(  \vec{0},2\varepsilon_{1}\right)  =\Gamma
\Delta$ where $\Delta$ is a subgroup of $G_{\infty}\left(  \vec{0}%
,2\varepsilon_{1}\right)  $ acting freely on $B^{3}\left(  \varepsilon
_{1}\right)  $. Then $B^{3}\left(  \varepsilon_{1}\right)  /\Delta$ is a
manifold (where we think of $G_{\infty}$ acting on the left) and $\left.
B^{3}\left(  \varepsilon_{1}\right)  \right/  G_{\infty}\left(  \vec
{0},2\varepsilon_{1}\right)  $ is isometric to $\left[  B^{3}\left(
\varepsilon\right)  /\Delta\right]  /\Gamma.$ For the rest of this subsection
we assume that $X_{\infty}$ is the GH limit of an essential ANSC sequence to
ensure that $G_{\infty}^{0}$ is abelian. We proceed with our construction of
$\Delta.$

Since Lie$\left(  G_{\infty}^{0}\right)  $ is abelian, there is a Lie algebra
decomposition Lie$\left(  \Gamma^{0}\right)  \oplus\mathfrak{h}^{\prime}$. Let
$k=\dim($Lie$(G_{\infty}))$. If we take $\varepsilon_{2}>0$ small so that
$\exp:$ Lie$\left(  \Gamma^{0}\right)  \cap B^{k}\left(  \varepsilon
_{2}\right)  \rightarrow G_{\infty}^{0}$ is injective, then we can define the
local group $\Delta=\left\{  \exp v:v\in\mathfrak{h}^{\prime}\cap B^{k}\left(
\varepsilon_{2}\right)  \right\}  .$ We note that $\Delta$ is a sub-local
group of $G_{\infty}^{0}$ and that for all $\delta\in\Delta-\left\{
id\right\}  ,$ $\delta\left(  \vec{0}\right)  \neq\vec{0}.$ Since Lie$\left(
\Gamma^{0}\right)  $ is abelian, $\left(  \exp v\right)  \left(  \exp
v^{\prime}\right)  =\exp\left(  v+v^{\prime}\right)  $ for $v\in\,$Lie$\left(
\Gamma^{0}\right)  $ and $v^{\prime}\in\mathfrak{h}^{\prime}$ so $\Gamma
^{0}\Delta$ generates a neighborhood $G_{\infty}^{0}\left(  \vec
{0},2\varepsilon_{0}\right)  $ of the identity for some small $\varepsilon
_{0}\in(0,1/4].$ We have the following properties of $\Delta$.

\begin{lemma}
There exists $\varepsilon_{0}\in(0,1/4]$ such that

\begin{enumerate}
\item $\Gamma^{0}\Delta=\Delta\Gamma^{0}=G_{\infty}^{0}\left(  \vec
{0},2\varepsilon_{0}\right)  ,$

\item $\Delta\cap\Gamma=\left\{  id\right\}  ,$ and

\item $\Delta$ acts freely on $B^{3}\left(  \varepsilon_{0}\right)  .$
\end{enumerate}
\end{lemma}

\begin{proof}
1 and 2 follow from the discussion above. To prove 3, suppose for every
$\varepsilon_{1}>0$, there exists $\delta\in\Delta$ and $p\in B^{3}\left(
\varepsilon_{1}\right)  $ such that $\delta\left(  p\right)  =p.$ Then we can
take a sequence $v_{i}\in\mathfrak{h}^{\prime}$ and $p_{i}\in B^{3}\left(
1/i\right)  $ such that $\left(  \exp v_{i}\right)  \left[  p_{i}\right]
=p_{i}.$ There exists a subsequence $i_{j}$ and $n_{j}\in\mathbb{N}$ such that
$n_{j}v_{i_{j}}\rightarrow v_{\infty}\neq0.$ Since $p_{i}\rightarrow\vec{0}$
and $\exp\left(  n_{i}v_{i_{j}}\right)  \left[  p_{i_{j}}\right]  =\left(
\exp v_{i_{j}}\right)  ^{n_{j}}\left[  p_{i_{j}}\right]  =p_{i_{j}}$ we have
$\left(  \exp v_{\infty}\right)  \left[  \vec{0}\right]  =\vec{0},$ so
$v_{\infty}\notin\mathfrak{h}^{\prime}.$ However, since $\mathfrak{h}^{\prime
}$ is closed, $v_{\infty}\in\mathfrak{h}^{\prime}$, a contradiction.
\end{proof}

Now we can decompose the entire action as follows.

\begin{lemma}
There exist $\varepsilon_{0}\in(0,1/4]$ and a manifold with Riemannian metrics
$\left(  V^{m},h\left(  t\right)  \right)  $ such that $\Gamma$ acts on
$\left(  V^{m},h\left(  t\right)  \right)  $ isometrically and $\left.
\left(  V^{m},h\left(  t\right)  \right)  \right/  \Gamma=\left.  \left(
B^{3}\left(  \varepsilon_{0}\right)  ,\tilde{g}_{\infty}\left(  t\right)
\right)  \right/  G_{\infty}\left(  \vec{0},2\varepsilon_{0}\right)  .$
\end{lemma}

\begin{proof}
By construction $G_{i}\left(  \vec{0},2\varepsilon_{0}\right)  $ preserves the
metrics $\tilde{g}_{i}\left(  t\right)  $ for all $t$ and hence $G_{\infty
}\left(  \vec{0},2\varepsilon_{0}\right)  $ preserves the metrics $\tilde
{g}_{\infty}\left(  t\right)  $ for all $t.$ Since $\Delta$ acts freely,
$V^{m}=\left.  B^{3}\left(  \varepsilon_{0}\right)  \right/  \Delta$ is a
manifold with quotient Riemannian metrics $h\left(  t\right)  $ induced by the
metrics $\tilde{g}_{\infty}\left(  t\right)  $. We claim that $\Delta$ is
normal in $G_{\infty}\left(  \vec{0},2\varepsilon_{0}\right)  ,$ this implies
that there is an isometric action $\Gamma\times V^{m}\rightarrow V^{m}.$ The
lemma follows from the fact that $\Gamma\Delta=G_{\infty}\left(  \vec
{0},2\varepsilon_{0}\right)  .$ It is clear that $\Gamma\subset O(m)$ where
$m=\dim V$ by looking at the derivatives of the action.

The claim can be proved from the fact that $G_{\infty}$ acts isometrically on
$D^{2}\left(  1/2\right)  \times(-1/2,1/2)$ with metric $^{2}\widetilde
{h}_{\infty}(t)+du^{2}$ (see \S\ref{sec 4.2 conse strong max}).
\end{proof}

We will call $\left(  V^{m},\Gamma,G_{\infty}\right)  $ the local model of
some neighborhood of $P_{\infty}$ in $X_{\infty}$.

\section{Classifying the limit space\label{Section 4 limit space}}

\subsection{General properties of GH limits of solutions}

Let $\left\{  \left(  M_{i}^{3},g_{i}\left(  t\right)  ,O_{i}\right)
\right\}  $ be a sequence with ANSC, essential collapsing split-like origins,
and unbounded diameter. By the Gromov compactness theorem there exists a
subsequence which converges in pointed GH distance to Alexandrov spaces
$(X_{\infty}(t),d_{\infty}(t),O_{\infty})$. We can use the fact that solutions
to the RF are uniformly bi-Lipschitz to each other to see that the metric
spaces $X_{\infty}(t)$ are topologically the same for all $t,$ and call that
space $X_{\infty}=X_{\infty}(0)$. Details are in \cite{Gl}. It follows from
\cite{Gl} that $\left.  \left(  B^{3}\left(  \varepsilon_{0}\right)
,\tilde{g}_{\infty}\left(  t\right)  \right)  \right/  G_{\infty}\left(
\vec{0},2\varepsilon_{0}\right)  $ in \S\ref{sec 3.3 local model} are
isometric to metric spaces $\left(  B_{d_{\infty}\left(  0\right)
}^{X_{\infty}}\left(  P_{\infty},\varepsilon_{0}\right)  ,d_{\infty
}(t)\right)  $ if $\varepsilon_{0}$ is taken small enough (independent of $t$).

$(X_{\infty},d_{\infty}(t),O_{\infty})$ has nonnegative curvature by the
following standard result: a GH limit of pointed Alexandrov spaces with
curvature $\geq$ k is itself a space of curvature $\geq$ k \cite[Proposition
10.7.1]{BBI}. We summarize in the following proposition.

\begin{proposition}
Let $\left\{  \left(  M_{i}^{3},g_{i}\left(  t\right)  ,O_{i}\right)
\right\}  $ be a sequence with ANSC, essential collapsing split-like origins,
and unbounded diameter. Then there is a subsequence which converges to an
Alexandrov spaces $\left(  X_{\infty},d_{\infty}\left(  t\right)  ,O_{\infty
}\right)  $ with nonnegative curvature whose dimension is $1$ or $2.$
\end{proposition}

\subsection{Consequences of the strong maximum principle for RF on local
covering geometries\label{sec 4.2 conse strong max}}

The ANSC condition and Hamilton's strong maximum principle for systems\ will
restrict the limit local covering geometries of the sequence. When combined
with low dimension of the sequence $\left\{  \left(  M_{i}^{3},g_{i}\left(
t\right)  ,O_{i}\right)  \right\}  $, it restricts how the local Lie groups of
isometries act.

Assume that all sectional curvatures $K_{g_{i}(0)}\leq1$ on $M_{i}^{3}$.
Choose $P_{i}$ as in \S\ref{sec 3.2 limit metric} to be $O_{i}$ and let
$\widetilde{g}_{O_{\infty}}(t)$ denote the limit metric. Then the\ essential
collapsing split-like origins imply that at origin $O_{\infty}$ and $t=0,$
$\lambda_{1}\left(  Rm_{\widetilde{g}_{O_{\infty}}(0)}(O_{\infty})\right)
=\lambda_{2}\left(  Rm_{\widetilde{g}_{O_{\infty}}(0)}(O_{\infty})\right)  =0$
and $0<\lambda_{3}\left(  Rm_{\widetilde{g}_{O_{\infty}}(0)}(O_{\infty
})\right)  $. By the strong maximum principle, the metrics $\widetilde
{g}_{O_{\infty}}\left(  t\right)  $ locally split as the product of a surface
metric with $\mathbb{R}.$ Moreover, the image of $Rm\left[  \widetilde
{g}_{O_{\infty}}\left(  t\right)  \right]  $ in $\wedge^{2}B^{3}(1)$ is
1-dimensional, independent of time, and invariant under parallel translation
(see \cite[Theorem 8.3]{Hamilton 4d PCO}). There is a unit 1-form which is
parallel and independent of time spanning the null space of\label{change 2}
$Rc\left[  \widetilde{g}_{O_{\infty}}\left(  t\right)  \right]  $ and
perpendicular to the image of $Rm\left[  \widetilde{g}_{O_{\infty}}\left(
t\right)  \right]  .$ By the deRham theorem, since $B^{3}\left(  1\right)  $
is contractible, for all $t\in\left(  \alpha,\omega\right)  $ there exists a
solution $^{2}\widetilde{h}_{O_{\infty}}(t)$ to the RF on $D^{2}\left(
1/2\right)  $ such that
\[
\iota_{_{O_{\infty}}}:\left(  D^{2}\left(  1/2\right)  \times(-1/2,1/2),^{2}%
\widetilde{h}_{O_{\infty}}(t)+du^{2},(\vec{0},0)\right)  \hookrightarrow
\left(  B^{3}(1),\widetilde{g}_{O_{\infty}}(t),\vec{0}\right)
\]
is an isometric embedding of the product of an evolving surface with an
interval and $^{2}\widetilde{h}_{O_{\infty}}(0)$ is a metric in normal
coordinates on the disk $D^{2}\left(  1/2\right)  .$

\label{following paragraph is mostly new}The local covering geometry of any
point $P_{\infty}\in X_{\infty}$ also splits as the product of a surface and
an interval. One way to see this is as follows. Let $P_{i}\in M_{i}$ be a
sequence of points such that $\left(  M_{i},d_{g_{i}\left(  t\right)  }%
,P_{i}\right)  $ converges to $\left(  X_{\infty},d_{\infty}\left(  t\right)
,P_{\infty}\right)  .$ Let $\gamma_{i}:\left[  0,L_{i}\right]  \rightarrow
M_{i}$ be a minimal geodesic joining $O_{i}$ to $P_{i}$ with respect to the
metric $g_{i}\left(  0\right)  .$ $\gamma_{i}$ extends (uniquely) to a
geodesic $\bar{\gamma}_{i}:\left[  -1,L_{i}+1\right]  \rightarrow M_{i}.$
Consider the geodesic tube $T_{i}:B^{2}\left(  1\right)  \times\left[
-1,L_{i}+1\right]  \rightarrow M_{i}$ corresponding to $\bar{\gamma}_{i}$ with
the pulled-back metrics $\bar{g}_{i}\left(  t\right)  =T_{i}^{\ast}%
g_{i}\left(  t\right)  ,$ which are solutions to the RF. By passing to a
subsequence, we obtain a limit solution $\left(  B^{2}\left(  1\right)
\times\left(  -1,L_{\infty}+1\right)  ,\bar{g}_{\infty}\left(  t\right)
\right)  $ with nonnegative sectional curvature, where $L_{\infty}%
=\lim_{i\rightarrow\infty}L_{i}=d_{\infty}\left(  0\right)  \left[  O_{\infty
},P_{\infty}\right]  $. Since the origins are split-like, $\lambda_{1}\left(
Rm_{\overline{g}_{\infty}\left(  t\right)  }\right)  =0$ at $\left(  \vec
{0},0\right)  $ for all $t\in\left(  \alpha,\omega\right)  ,$ where $\vec{0}$
is the origin in $B^{2}\left(  1\right)  .$ By the strong maximum principle,
$\lambda_{1}\left(  Rm_{\overline{g}_{\infty}\left(  t\right)  }\right)  =0$
at $\left(  \vec{0},L_{\infty}\right)  ;$ actually $\lambda_{1}\left(
Rm_{\overline{g}_{\infty}\left(  t\right)  }\right)  =0$ for every point in
$B^{2}\left(  1\right)  \times\left(  -1,L_{\infty}+1\right)  $. This implies
that $\lim_{i\rightarrow\infty}\lambda_{1}\left(  Rm_{g_{i}}(P_{i})\right)
=0$ for all $t$ and hence $\lambda_{1}\left(  Rm_{\widetilde{g}_{P_{\infty}}%
}\left(  \tilde{P}_{\infty}\right)  \right)  =0$ where $\tilde{P}_{\infty}$ is
the origin of balls $B^{3}\left(  1\right)  $ and $\widetilde{g}_{P_{\infty}%
}\left(  t\right)  $ is the limit metric coming from the local covering
geometry construction around the $P_{i}$. From this we conclude that the local
covering geometry $\widetilde{g}_{P_{\infty}}\left(  t\right)  $ of
$P_{\infty}\in X_{\infty}$ has a zero curvature and hence splits as a product
of a surface and an interval, similarly to how it splits at $O_{\infty}.$

Another way to prove $\lambda_{1}\left(  Rm_{\widetilde{g}_{P_{\infty}}\left(
t\right)  }\left(  \tilde{P}_{\infty}\right)  \right)  =0$ is to assume, by
contradiction, that there is a subsequence such that the $P_{i}$'s are
bump-like. Then there exists a uniform injectivity radius estimate at $P_{i}.$
Since $d_{g_{i}\left(  0\right)  }\left(  P_{i},O_{i}\right)  $ is uniformly
bounded, this implies a uniform injectivity radius estimate at $O_{i},$ which
is a contradiction.

\subsection{Killing vector fields on surfaces\label{B4bbb}}

We shall show that a local surface with a Killing vector field is locally a
warped product. This clearly must be a classical fact\label{find reference 1}
but since we have not found a reference we include a sketch of the proof.

Let $\left(  \Sigma^{2},h\right)  $ be a Riemannian surface (not necessarily
complete) with a Killing vector field $K$. Let $J:T\Sigma^{2}\rightarrow
T\Sigma^{2}$ be the complex structure, that is, rotation by $90^{\circ}$ in
the counterclockwise direction. Let $x\in\Sigma^{2}$ and define a smooth unit
speed path $\gamma:\left(  r_{0}-\varepsilon,r_{0}+\varepsilon\right)
\rightarrow\Sigma^{2}$ by $\dot{\gamma}\left(  r\right)  =\frac{J\left(
K\right)  }{\left|  J\left(  K\right)  \right|  }\left(  \gamma\left(
r\right)  \right)  ,$ $\gamma\left(  r_{0}\right)  =x$. Define also a
1-parameter family of smooth paths $\beta_{r}:\left(  \theta_{0}%
-\varepsilon,\theta_{0}+\varepsilon\right)  \rightarrow\Sigma^{2}$ by
$\dot{\beta}_{r}\left(  \theta\right)  =K\left(  \beta_{r}\left(
\theta\right)  \right)  ,$ $\beta_{r}\left(  \theta_{0}\right)  =\gamma\left(
r\right)  ,$ and a 1-parameter family of smooth unit speed paths
$\gamma_{\theta}:\left(  r_{0}-\varepsilon,r_{0}+\varepsilon\right)
\rightarrow\Sigma^{2}$ by $\dot{\gamma}_{\theta}\left(  r\right)
=\frac{J\left(  K\right)  }{\left|  J\left(  K\right)  \right|  }\left(
\gamma_{\theta}\left(  r\right)  \right)  ,$ $\gamma_{\theta}\left(
r_{0}\right)  =\beta_{r_{0}}\left(  \theta\right)  $. Note that $\gamma
_{\theta_{0}}=\gamma.$

\begin{lemma}%
\[
\beta_{r}\left(  \theta\right)  =\gamma_{\theta}\left(  r\right)  .
\]
\end{lemma}

\begin{proof}
This follows from%
\[
\left[  K,\frac{J\left(  K\right)  }{\left|  J\left(  K\right)  \right|
}\right]  =\frac{1}{\left|  J\left(  K\right)  \right|  }\left[  K,J\left(
K\right)  \right]  -\frac{1}{2}\frac{1}{\left|  J\left(  K\right)  \right|
^{3}}K\left|  J\left(  K\right)  \right|  ^{2}J\left(  K\right)  =0
\]
since $\left[  K,J\left(  K\right)  \right]  =0$ and $K\left|  J\left(
K\right)  \right|  ^{2}=0.$
\end{proof}

Hence $(r,\theta)$ defines local coordinates on a neighborhood of $x$. Define
the function%
\[
f\left(  r\right)  \doteqdot\left|  K\right|  \left(  \beta_{r}\right)
\]
where we are using the fact that $\left|  K\right|  $ is constant on
$\beta_{r}$. The metric is given by%
\begin{align*}
h  &  =\left\langle \frac{J\left(  K\right)  }{\left|  J\left(  K\right)
\right|  },\frac{J\left(  K\right)  }{\left|  J\left(  K\right)  \right|
}\right\rangle dr^{2}+2\left\langle K,\frac{J\left(  K\right)  }{\left|
J\left(  K\right)  \right|  }\right\rangle drd\theta+\left\langle
K,K\right\rangle d\theta^{2}\\
&  =dr^{2}+f\left(  r\right)  ^{2}d\theta^{2}.
\end{align*}
We have proved:

\begin{lemma}
\label{lem coord on surf}(i) Given $x\in\Sigma^{2}$, there exists a
neighborhood $U$ of $x$ and local coordinates $r$ and $\theta$ on $U$ such
that%
\begin{align*}
\frac{\partial}{\partial r}  &  =\frac{J\left(  K\right)  }{\left|  J\left(
K\right)  \right|  },\text{ \ \ }\frac{\partial}{\partial\theta}=K,\\
h  &  =dr^{2}+f\left(  r\right)  ^{2}d\theta^{2}.
\end{align*}
(ii) Let $(r_{i},\theta_{i}),i=1,2$ be two coordinate systems of some
neighborhood $U$ of $x$ and $f_{i}(r_{i})$ be the functions defining the
metric in (i), then
\begin{align*}
r_{2}  &  =r_{1}+r_{0},\text{ \ \ }\theta_{2}=\theta_{1}+\theta_{0},\\
f_{2}(r_{2})  &  =f_{1}(r_{2}-r_{0})
\end{align*}
where $r_{0}=r_{2}(x)-r_{1}(x)$ and\ $\theta_{0}=\theta_{2}(x)-\theta_{1}(x)$
are constants.
\end{lemma}

Next we prove a uniqueness theorem about Killing vector fields.

\begin{lemma}
\label{Lem Abelian act}Let $\left(  \Sigma^{2},h\right)  $ be a Riemannian
surface with nonzero curvature everywhere. Suppose $K_{1}$ and $K_{2}$ are two
Killing vector fields satisfying $\left[  K_{1},K_{2}\right]  \equiv0$. Then
$K_{1}$ and $K_{2}$ are linearly dependent.
\end{lemma}

\begin{proof}
Using $\left[  K_{1},K_{2}\right]  \equiv0$ and the Killing vector field
equation for $K_{i}$
\begin{equation}
\left\langle \nabla_{X}K_{i},Y\right\rangle +\left\langle \nabla_{Y}%
K_{i},X\right\rangle =0 \label{Killing vf eqn for K_i X Y}%
\end{equation}
we have $\left\langle \nabla_{K_{2}}K_{1},K_{2}\right\rangle =0$ and
$\left\langle \nabla_{K_{2}}K_{1},K_{1}\right\rangle =\left\langle
\nabla_{K_{1}}K_{2},K_{1}\right\rangle =0$. If $K_{1}$ and $K_{2}$ are
linearly independent, then $\nabla_{K_{2}}K_{1}\equiv\nabla_{K_{1}}K_{2}%
\equiv0.$ Also from $\left\langle \nabla_{K_{1}}K_{1},K_{2}\right\rangle
=-\left\langle \nabla_{K_{2}}K_{1},K_{1}\right\rangle =0$ and $\left\langle
\nabla_{K_{1}}K_{1},K_{1}\right\rangle =0$ we get $\nabla_{K_{1}}K_{1}\equiv
0$. Similarly we have $\nabla_{K_{2}}K_{2}\equiv0$. Hence the metric is flat,
a contradiction.
\end{proof}

\subsection{A canonical form for actions on a surface $\times$ $\mathbb{R}%
$\label{eee}}

Let $(\Sigma^{2},h)$ be a Riemannian surface with positive curvature.

\begin{lemma}
\label{structure of group action}Given an $\mathbb{R}_{\text{loc}}$ local
group action of local isometries on $\Sigma^{2}\times\mathbb{R}_{\text{loc}}$
with the product metric, there exist coordinates $(r,\theta)$ on $\Sigma^{2}$
such that we can write the action of $\tau\in\mathbb{R}_{\text{loc}}$ as%
\begin{equation}
\tau:\left(  r,\theta,u\right)  \mapsto\left(  r,\theta+a\tau,u+b\tau\right)
\label{S^1 x R local action form}%
\end{equation}
for some constant $a,b\in\mathbb{R}$.
\end{lemma}

\begin{proof}
Let $\left(  x,y\right)  $ be coordinates on $\Sigma^{2}$. Denote the image of
$\left(  x,y,u\right)  $ under action of $\tau$ by $(x_{\tau},y_{\tau}%
,u_{\tau}).$ Since $(\Sigma^{2},h)$ has positive curvature and $\tau$ is a
local isometry, the tangential map $d\tau$ must be of the block diagonal form
\[
T_{\left(  x,y\right)  }\Sigma^{2}\times T_{u}\mathbb{R}_{\text{loc}%
}\rightarrow T_{\left(  x_{\tau},y_{\tau}\right)  }\Sigma^{2}\times
T_{u_{\tau}}\mathbb{R}_{\text{loc}}.
\]
This implies that the functions $x_{\tau}(x,y,u)$ and $y_{\tau}(x,y,u)$ are
independent of $u$ and the function $u_{\tau}(x,y,u)$ is independent of
$(x,y)$. We can write $x_{\tau}(x,y,u)=x_{\tau}(x,y)$, $y_{\tau}%
(x,y,u)=y_{\tau}(x,y),$ and $u_{\tau}(x,y,u)=u_{\tau}(u)$.

The assumption that the action is locally isometric on $\Sigma^{2}%
\times\mathbb{R}_{\text{loc}}$ implies that the action $\tau:\mathbb{R}%
_{\text{loc}}\rightarrow\mathbb{R}_{\text{loc}}$ with $\tau:u\rightarrow$
$u_{\tau}(u)$ is an isometry and the action $\tau:\Sigma^{2}\rightarrow
\Sigma^{2}$ with $\tau:(x,y)\rightarrow(x_{\tau}(x,y),y_{\tau}(x,y))$ is an
isometry. Hence $u_{\tau}(u)=\pm u+b\tau$. If action $\tau:\Sigma
^{2}\rightarrow\Sigma^{2}$ is nontrivial, then by Lemma \ref{lem coord on
surf} there are coordinates $(r,\theta)$ on $\Sigma^{2}$ such that
$\tau:\left(  r,\theta\right)  \mapsto\left(  r,\pm\theta+a\tau\right)  $ for
some $a\neq0$. If the action $\tau:\Sigma^{2}\rightarrow\Sigma^{2}$ is
trivial, then any coordinates on $\Sigma^{2}$ will be fine with $a=0.$
\end{proof}

\subsection{Listing of the possible local model data $(V^{m},\Gamma,G_{\infty
})$ \label{B1c}}

\textbf{\label{priority for working on 1}}In this section we assume that
$X_{\infty}$ is the GH limit of an essential ANSC sequence to ensure that
$G_{\infty}^{0}$ is abelian. We shall classify all possible local model data
$(V^{m},\Gamma,G_{\infty})$ for a neighborhood of $P_{\infty}\in X_{\infty}$.

The Hausdorff dimension of the limit $X_{\infty}$ is the same everywhere
\cite[Theorem 10.6.1]{BBI} and must be an integer \cite[Theorem 10.8.2]{BBI}.
Hence $\dim X_{\infty}$ is either $1$ or $2$. The dimension cannot be $0$
since the diameters are unbounded and the dimension cannot be $3$ since the
origins are collapsing.

Recall the following elementary fact \cite[Theorem 2.2.1]{GB}.

\begin{lemma}
The only discrete subgroups of $O\left(  2\right)  $ are $\mathbb{Z}_{p}$ and
$\mathbb{D}_{2p}=\mathbb{Z}_{p}\rtimes\mathbb{Z}_{2},$ the dihedral group of
order $2p,$ and the only one dimensional subgroups are $SO\left(  2\right)  $
and $O\left(  2\right)  .$
\end{lemma}

We can now enumerate the possible models; we will use the notations from
\S\ref{sec 3.3 local model}. It will be important to check that $\Delta$ is
normal in $G_{\infty}^{0},$ which we will check individually for each case.

\textbf{Case 0}. ($m=0$) This case cannot occur because of the unbounded
diameters assumption.

\textbf{Case 1}. ($m=1$) It follows from the ANSC assumption and \S\S\ref{sec
4.2 conse strong max} and \ref{eee} that there is coordinate $\left(
r,\theta,u\right)  $ on $D^{2}\left(  1/2\right)  \times(-1/2,1/2)$ with
metric $\tilde{g}_{\infty}\left(  t\right)  =dr^{2}+f^{2}(r,t)d\theta
^{2}+du^{2},$ such that $\Delta\ni(\tau_{1},\tau_{2})$ acts by $\left(
r,\theta,u\right)  \rightarrow(r,\theta+\tau_{1},u+\tau_{2})$. $\Gamma$ can be either

(\textbf{1i}) $\left\{  0\right\}  ,$ or

(\textbf{1ii}) $\mathbb{Z}_{2}$ acting by $-1:\left(  r,\theta,u\right)
\rightarrow\left(  r,-\theta,-u\right)  $, or

(\textbf{1iii}) $-1:\left(  r,\theta,u\right)  \rightarrow\left(
-r,\theta,-u\right)  $ when $f(r)=f(-r),$\ or

(\textbf{1iv}) $-1:\left(  r,\theta,u\right)  \rightarrow\left(
-r,-\theta,u\right)  $ when $f(r)=f(-r)$.

In all cases $\Delta$ is normal in $G_{\infty}\left(  \vec{0},2\varepsilon
_{0}\right)  .$ This implies that $P_{\infty}$ has a neighborhood homeomorphic
to $\left(  a,b\right)  ,$ for some $a<0<b$ in case (1i) and (1ii), and
$P_{\infty}$ has a neighborhood homeomorphic to $[0,b),$ for some $0<b$ in
case (1iii) and (1iv). $P_{\infty}$ corresponds to $0$. $G_{\infty}%
^{0}=\mathbb{R}_{loc}^{2}$ and $G_{\infty}=\mathbb{R}_{loc}^{2}$ in case (1i)
and $G_{\infty}=\mathbb{R}_{loc}^{2}\rtimes\mathbb{Z}_{2}$ in case (1ii),
(1iii) and (1iv). We will rule out case (1iii) and (1iv) in the proof of
Theorem \ref{Thm 1D virt mod}.

\textbf{Case 2a}. ($m=2$, $\dim X_{\infty}=1$) Then $\dim\Gamma^{0}=1$ and the
action generates a Killing vector field on $D^{2}\left(  1/2\right)  $ which
vanishes at origin. It follows from the ANSC assumption and \S\S\ref{sec 4.2
conse strong max} and \ref{eee} that there is coordinate $\left(
r,\theta,u\right)  $ on $D^{2}\left(  1/2\right)  \times(-1/2,1/2)$ with
metric $\tilde{g}_{\infty}\left(  t\right)  =dr^{2}+f^{2}(r,t)d\theta
^{2}+du^{2}$ and $f(0,t)=0,$ such that $\Gamma^{0}\ni\tau$ acts by $\left(
r,\theta,u\right)  \rightarrow(r,\theta+\tau,u)$. $\Gamma$ can be either

(\textbf{2ai}) $\Gamma=SO(2)=\Gamma^{0}$ or

(\textbf{2aii}) $\Gamma=O(2)=\Gamma^{0}\rtimes\mathbb{Z}_{2}$ acting by
$-1:\left(  r,\theta,u\right)  \rightarrow\left(  r,-\theta,-u\right)  $ for
$-1\in\mathbb{Z}_{2}.$

In both cases $\Delta$ is normal in $G_{\infty}\left(  \vec{0},2\varepsilon
_{0}\right)  $ and acts by $\Delta\ni\tau:(r,\theta,u)\rightarrow
(r,\theta+c_{1}\tau,u+c_{2}\tau)$ with two constants $c_{1}$ and $c_{2}\neq0.$
This implies that $P_{\infty}$ has a neighborhood homeomorphic to $[0,b),$ for
some $b>0$ and $P_{\infty}$ corresponds to $0$ in both cases. $G_{\infty}%
^{0}=\mathbb{R}_{loc}^{2}$ in both case. $G_{\infty}=\mathbb{R}_{loc}\times
SO(2)$ in case (2ai) and $G_{\infty}=\mathbb{R}_{loc}\times O(2)$ in case (2aii).

\textbf{Case 2b}. ($m=2$, $\dim X_{\infty}=2$) Then $\Gamma$ is discrete;
since $\Gamma^{0}=\{id\}$ so $\Delta$ is normal in $G_{\infty}\left(  \vec
{0},2\varepsilon_{0}\right)  $ and $\Gamma\subset O(2)$. By the lemma, either

(\textbf{2bi}) $\Gamma=\mathbb{Z}_{p}$ for some $p\in\mathbb{N}$ or

(\textbf{2bii}) $\Gamma=\mathbb{D}_{2p}.$

This implies that $P_{\infty}$ has a neighborhood homeomorphic to (2bi)
$D^{2}/\mathbb{Z}_{p}$ or (2bii) $D^{2}/\mathbb{D}_{2p}$ where $D^{2}$ is a
two dimensional disk and $\mathbb{Z}_{2}<\mathbb{D}_{2p}$ acts by reflection.
$G_{\infty}^{0}=\mathbb{R}_{loc}.$

\textbf{Case 3}. ($m=3$) This case cannot occur. By Lemma \ref{group
infinite}, $\dim\Delta\geq1$, so $m=\dim V=3-\dim\Delta\leq2.$ This is a contradiction.

We summarize the models in the following proposition$.$

\begin{proposition}
\label{local model}The following tables give the complete list of the
possibilities for local models of collapse in our situation (where
$p\in\mathbb{N}$):%
\[%
\begin{tabular}
[c]{|l|l|l|l|}\hline
$m=1$ & $\Gamma$ & ${\Large \mathstrut}G_{\infty}^{0}$ & $X_{\infty}$
(locally)\\\hline
(1i) & $\left\{  0\right\}  $ & ${\Large \mathstrut}\mathbb{R}_{loc}^{2}$ &
$(a,b)$\\\hline
(1ii) & $\mathbb{Z}_{2}$ & ${\Large \mathstrut}\mathbb{R}_{loc}^{2}$ &
$(a,b)$\\\hline
(1iii) & $\mathbb{Z}_{2}$ & ${\Large \mathstrut}\mathbb{R}_{loc}^{2}$ &
$[0,b)$\\\hline
(1iv) & $\mathbb{Z}_{2}$ & ${\Large \mathstrut}\mathbb{R}_{loc}^{2}$ &
$[0,b)$\\\hline
\end{tabular}
\]
and
\[%
\begin{tabular}
[b]{|l|l|l|l|}\hline
$m=2$ & $\Gamma$ & ${\Large \mathstrut}G_{\infty}^{0}$ & $X_{\infty}$
(locally)\\\hline
(2ai) & $SO\left(  2\right)  $ & ${\Large \mathstrut}\mathbb{R}_{loc}%
^{1}\times SO(2)$ & $[0,b)$\\\hline
(2aii) & $O(2)$ & ${\Large \mathstrut}\mathbb{R}_{loc}^{1}\times SO(2)$ &
$[0,b)$\\\hline
(2bi) & $\mathbb{Z}_{p}$ & ${\Large \mathstrut}\mathbb{R}_{loc}$ &
$D^{2}/\mathbb{Z}_{p}$\\\hline
(2bii) & $\mathbb{D}_{2p}$ & ${\Large \mathstrut}\mathbb{R}_{loc}$ &
$D^{2}/\mathbb{D}_{2p}$\\\hline
\end{tabular}
.
\]
\end{proposition}

\section{Constructing the 2-dimensional virtual limit\label{sec virtual}}

In this section $(X_{\infty},h_{\infty}(t),O_{\infty})$ is the limit of a
sequence $\left\{  \left(  M_{i}^{3},g_{i}\left(  t\right)  ,O_{i}\right)
\right\}  $\ with ANSC, essential collapsing split-like origins, and unbounded diameters.

\subsection{The 2-dimensional limit orbifold and its virtual limit \label{jjj}}

First we recall the following definition of orbifolds ($V$- manifolds) due to
I. Satake (1956) and W. Thurston.

\begin{definition}
An $n$-dimensional orbifold $\mathcal{X}$ is a Hausdorff space $X$ together
with a collection of pairs of open sets and finite groups $\left\{  \left(
U_{i},\Gamma_{i}\right)  \right\}  $ such that

\begin{enumerate}
\item $\left\{  U_{i}\right\}  $ is closed under finite intersections;

\item  For each $i,$ there is an open subset $\tilde{U}_{i}\subset
\mathbb{R}^{n}$ such that $\Gamma_{i}$ acts on $\tilde{U}_{i}$ and
$U_{i}\approx\tilde{U}_{i}/\Gamma_{i}$;

\item  Whenever $U_{i}\subset U_{j}$ there is an injective homomorphism
$f_{ij}:\Gamma_{i}\rightarrow\Gamma_{j}$ and an embedding $\tilde{\phi}%
_{ij}:\tilde{U}_{i}\hookrightarrow\tilde{U}_{j}$ such that
\[
\tilde{\phi}_{ij}\left(  \gamma x\right)  =f_{ij}\left(  \gamma\right)
\tilde{\phi}_{ij}\left(  x\right)
\]
for all $\gamma\in\Gamma_{i}$ and such that the appropriate diagram commutes.
An orbifold is smooth if $\tilde{\phi}_{ij}$ are smooth maps.
\end{enumerate}
\end{definition}

Recall that a length space such that every point has a neighborhood isometric
to a region with a Riemannian metric is a smooth Riemannian manifold (see, for
instance, \cite[\S5.1]{BBI}). We give an analogous result for a smooth orbifold.

\begin{proposition}
\label{orbifold characterization}If $X$ is a length space such that every
point $p\in X$ has a neighborhood $U$ isometric to $\tilde{U}_{p}/\Gamma_{p}$
where $\tilde{U}_{p}\subset\mathbb{R}^{n}$ (with a smooth Riemannian metric)
is a simply connected neighborhood of the origin $\vec{0}$ and $\Gamma_{p}$ is
a finite group acting effectively on $\tilde{U}_{p}$ such that

\begin{enumerate}
\item  the isometries $\tilde{U}_{p}/\Gamma_{p}\rightarrow U_{p}$ take
$\vec{0}$ to $p$ and either

\item[2a.] the only fixed point of the action of $\Gamma_{p}$ is $\vec{0}$ (if
there is any), or

\item[2b.] $\Gamma_{p}=\Gamma_{p}^{\prime}\rtimes\mathbb{Z}_{2}$ where the
only fixed point of the action $\Gamma_{p}^{\prime}$ is $\vec{0}$,
$\mathbb{Z}_{2}$ acts by reflection with respect to a hyperplane and
$\Gamma_{p}^{\prime}<\Gamma_{p}^{\prime}\rtimes\mathbb{Z}_{2}$ is a normal subgroup.
\end{enumerate}

Then $X,$ together with the open sets $\left\{  \left(  U_{p},\Gamma
_{p}\right)  \right\}  $ and all intersections of these open sets paired with
the trivial group or $\mathbb{Z}_{2},$ is a smooth Riemannian orbifold.
\end{proposition}

\begin{proof}
Since the only fixed points of the group $\Gamma_{p}^{\prime}$ actions are at
the origin in $\tilde{U}_{p},$ whenever there is an inclusion
$U\hookrightarrow U^{\prime},$ the group for $U$ is trivial or $\mathbb{Z}%
_{2}$. Furthermore, the lift $\widetilde{U}\rightarrow\tilde{U}^{\prime}$ is
well defined. Since it is an isometry into the image, the map is smooth by
Myers-Steenrod (see, for instance, \cite[Theorem 9.1]{Pet}), and hence the
orbifold is smooth.
\end{proof}

\begin{lemma}
\label{1 co-D collapse}The two dimensional limit $\left(  X_{\infty}%
,d_{\infty}(t)\right)  $ is a smooth orbifold with singularities of types
$D^{2}/\mathbb{Z}_{p}$ or $D^{2}/\mathbb{D}_{2p}$, where $D^{2}\subset
\mathbb{R}^{2}$ is a disk.
\end{lemma}

\begin{proof}
By Proposition \ref{local model} every point has a neighborhood isometric to
$V/\Gamma$ where $V$ is a neighborhood of $\vec{0}$ in $\mathbb{R}^{2}$ and
$\Gamma$ is $\mathbb{Z}_{p}$ or $\mathbb{D}_{2p}$. By Proposition
\ref{orbifold characterization} this is a smooth orbifold with the stated
singularity types.
\end{proof}

Since $X_{\infty}$ is an orbifold and the metric $d_{\infty}\left(  t\right)
$ comes from a Riemannian structure, we can write the limit as $\left(
X_{\infty},h_{\infty}\left(  t\right)  \right)  $ where $h_{\infty}\left(
t\right)  $ is the Riemannian metric on the orbifold. Next we will prove:

\begin{lemma}
\label{lem pos curv orb}The two dimensional limit $\left(  X_{\infty
},h_{\infty}(t)\right)  $ is a positively curved orbifold.
\end{lemma}

\begin{proof}
At each point $P_{\infty}\in X_{\infty}$ there is a neighborhood of
$P_{\infty}$ isometric to a finite quotient of $\left(  V,h_{V}\left(
t\right)  \right)  =\left.  \left(  D^{2}\left(  1/2\right)  \times
(-1/2,1/2),^{2}\widetilde{h}_{\infty}(t)+du^{2}\right)  \right/  \Delta,$ with
$\Delta=\mathbb{R}_{\text{loc}}$ acting freely by isometries. By Lemma
\ref{structure of group action} there are coordinates $\left(  r,\theta
\right)  $ on $D^{2}\left(  1/2\right)  $ such that the Killing vector field
of the action of $\Delta$ is $K=\left(  a\frac{\partial}{\partial\theta
},b\frac{\partial}{\partial u}\right)  $ and $^{2}\widetilde{h}_{\infty
}(0)=dr^{2}+f_{P_{\infty}}\left(  r\right)  ^{2}d\theta^{2}$ has positive
curvature. We consider two cases.

A. If $b\neq0$, by a simple calculation using O'Neill's formula for the
submersion
\[
D^{2}\left(  1/2\right)  \times(-1/2,1/2)\rightarrow\left.  \left(
D^{2}\left(  1/2\right)  \times(-1/2,1/2)\right)  \right/  \Delta
\]
we conclude that the curvature of $\left(  V,h_{V}\left(  0\right)  \right)  $
is strictly positive everywhere.

B. If $b=0$, then $\left(  V,h_{V}\left(  0\right)  \right)  \cong\left(
\left(  D^{2}\left(  1/2\right)  ,^{2}\widetilde{h}_{\infty}(0)\right)
/\Delta\right)  \times(-1/2,1/2)$ has zero curvature everywhere.

Combining A and B we see that if we can find a point $P_{\infty}\in X_{\infty
}$ such that $b\neq0$ then for every point in $X_{\infty}$ we have $b\neq0.$
Hence $\left(  X_{\infty},h_{\infty}(0)\right)  $ has positive curvature and
so does $\left(  X_{\infty},h_{\infty}(t)\right)  $. Otherwise, for every
point in $X_{\infty}$ we have $b=0$ and $\left(  X_{\infty},h_{\infty
}(0)\right)  $ is flat. Since the $\Delta$ action is free, $f_{P_{\infty}%
}(r)>0$ when $b=0$. In the next two lemmas we will show by contradiction that
$\left(  X_{\infty},h_{\infty}(0)\right)  $ cannot be flat. More precisely we
will construct a complete metric on $\left(  -\infty,\infty\right)  \times
S^{1}$ of positive curvature in Lemma \ref{lem piece metric}; such a metric on
$\left(  -\infty,\infty\right)  \times S^{1}$ cannot exist by the well-known
theorem of Gromoll-Meyer. This gives the required contradiction.
\end{proof}

First we prove a general property of the limit metric $\widetilde
{g}_{P_{\infty}}(t)$. Let $\left\{  \left(  N_{i}^{n},g_{i},O_{i}\right)
\right\}  $ be a sequence of pointed complete $n$-dimensional Riemannian
manifolds with bounded geometry. We assume that $\left(  N_{i}^{n},g_{i}%
,O_{i}\right)  $ converges to the metric space $\left(  X_{\infty},d_{\infty
},O_{\infty}\right)  $ in the pointed GH topology. Fix an $\epsilon>0$ and let
$P_{\infty,\mu}$ and $P_{\infty,\nu}$ be two points in $X_{\infty}$.

By the definition of GH-convergence we can pass to a subsequence such that
there exist maps $\varphi_{i}:\left(  N_{i}^{n},d_{g_{i}},O_{i}\right)
\rightarrow\left(  X_{\infty},d_{\infty},O_{\infty}\right)  $ and $\psi
_{i}:\left(  X_{\infty},d_{\infty},O_{\infty}\right)  \rightarrow\left(
N_{i}^{n},d_{g_{i}},O_{i}\right)  $ which are $1/i$-pointed GH approximations.
Let $P_{i,\mu}\doteqdot\psi_{i}\left(  P_{\infty,\mu}\right)  $ and $P_{i,\nu
}\doteqdot\psi_{i}\left(  P_{\infty,\nu}\right)  $. Fix orthonormal frames
$F_{i,\mu}$ of $T_{P_{i,\mu}}N_{i}^{n}$ and $F_{i,\nu}$ of $T_{P_{i,\nu}}%
N_{i}^{n}$, and define the exponential maps restricted to balls:
\begin{align*}
\exp_{P_{i,\mu}}  &  :B_{\mu}^{n}(1+\epsilon)\rightarrow N_{i}^{n},\\
\exp_{P_{i,\nu}}  &  :B_{\nu}^{n}(1+\epsilon)\rightarrow N_{i}^{n},
\end{align*}
where $B_{\mu}^{n}(1+\epsilon)$ is a ball of radius $1+\epsilon$ in
$\mathbb{R}^{n}$ and $\mu$ is used to indicate the dependence on the sequence
and the point in $X_{\infty}.$ Assume that the sectional curvatures $K_{g_{i}%
}\leq1$ on the balls $B_{g_{i}}\left(  P_{i,\mu},1+\epsilon\right)  $ and
$B_{g_{i}}\left(  P_{i,\nu},1+\epsilon\right)  $ so that $\exp_{P_{i,\mu}}$
and $\exp_{P_{i,\nu}}$ are local diffeomorphisms. We assume that the
pulled-back metrics $\widetilde{g}_{_{i,\mu}}\doteqdot\exp_{P_{i,\mu}}^{\ast
}g_{i}$ converge to $\widetilde{g}_{\infty,\mu}$ on $B_{\mu}^{n}\left(
1\right)  $ and $\widetilde{g}_{i,\nu}\doteqdot\exp_{P_{_{i,\nu}}}^{\ast}%
g_{i}$ converges to $\widetilde{g}_{\infty,\nu}$ on $B_{\nu}^{n}\left(
1\right)  $.

The next lemma shows how minimal geodesics allow us to identify parts of balls
in the tangent space to create overlap maps. Consider a manifold $\left(
M^{n},g\right)  $ with pullback metrics with local covers $\exp_{P_{\mu}%
}:\left(  B_{\mu}^{n}\left(  1\right)  ,\vec{0}_{\mu}\right)  \rightarrow
\left(  M^{n},P_{\mu}\right)  $ and $\exp_{P_{\nu}}:\left(  B_{\nu}^{n}\left(
1\right)  ,\vec{0}_{\nu}\right)  \rightarrow\left(  M^{n},P_{\nu}\right)  $
(say, assume $K_{g}\leq1$), and define $\tilde{g}_{\mu}=\exp_{P_{\mu}}^{\ast
}g$ and $\tilde{g}_{\nu}=\exp_{P_{\nu}}^{\ast}g.$

\begin{lemma}
\label{lem overlap from geodesic}Let $\beta$ be a unit speed minimal geodesic
joining $P_{\mu}$ and $P_{\nu}$ such that $\beta\subset\exp_{P_{\mu}}\left(
B_{\mu}^{n}\left(  1\right)  \right)  \cup$ $\exp_{P_{\nu}}\left(  B_{\nu}%
^{n}\left(  1\right)  \right)  .$ Let $Q$ be the midpoint of $\beta.$ Then for
$\delta<1-d_{g}\left(  P_{\mu},P_{\nu}\right)  /2$ there exist $Q_{\mu}\in
B_{\mu}^{n}\left(  1\right)  ,$ and $Q_{\nu}\in B_{\nu}^{n}\left(  1\right)  $
such that $\left(  B^{n}\left(  Q_{\mu},\delta\right)  ,\tilde{g}_{\mu
}\right)  $ is isometric to $\left(  B^{n}\left(  Q_{\nu},\delta\right)
,\tilde{g}_{\nu}\right)  $ (the balls are taken with respect to the given metric).
\end{lemma}

\begin{proof}
There is a unique lift $\widetilde{\beta}_{\mu}$ in $\left(  B_{\mu}%
^{n}\left(  1\right)  ,\tilde{g}_{\mu}\right)  $ of $\beta$ with
$\widetilde{\beta}_{\mu}(0)=\vec{0}_{\mu}\in B_{\mu}^{n}\left(  1\right)  $
and a unique lift $\widetilde{\beta}_{\nu}$ in $\left(  B_{\nu}^{n}\left(
1\right)  ,\tilde{g}_{\nu}\right)  $ of $-\beta$ (which is $\beta$ with the
time parameter reversed) with $\widetilde{\beta}_{\nu}(0)=\vec{0}_{\nu}\in
B_{\nu}^{n}\left(  1\right)  $. $Q$ is lifted to $Q_{\mu}$ in $B_{\mu}%
^{n}\left(  1\right)  $ by $\widetilde{\beta}_{\mu}$ and $Q$ is lifted to
$Q_{\nu}$ in $B_{\nu}^{n}\left(  1\right)  $ by $\widetilde{\beta}_{\nu}.$
Note that we have chosen $\delta$ such that $\exp_{P_{\mu}}\left(
B^{n}\left(  Q_{\mu},\delta\right)  \right)  $ and $\exp_{P_{\nu}}\left(
B^{n}\left(  Q_{\nu},\delta\right)  \right)  $ are contained in that
$\exp_{P_{\mu}}\left(  B_{\mu}^{n}\left(  1\right)  \right)  \cap\exp_{P_{\nu
}}\left(  B_{\nu}^{n}\left(  1\right)  \right)  .$ Since $\exp_{P_{\mu}}$ and
$\exp_{P_{\nu}}$ are local diffeomorphisms we can lift the exponential map
$\exp_{P_{\mu}}$ to a unique diffeomorphism $\iota:B^{n}\left(  Q_{\mu}%
,\delta\right)  \rightarrow B^{n}\left(  Q_{\nu},\delta\right)  $ which
satisfies $\iota\left(  Q_{\mu}\right)  =Q_{\nu}.$ Moreover, the metrics are
mapped isometrically.
\end{proof}

We can now apply this to our setting.

\begin{lemma}
\label{lem isom metr Fuk}If $d_{\infty}(P_{\infty,\mu},P_{\infty,\nu})<2$,
then there are neighborhoods $B^{n}(Q_{\infty,\mu},\delta)\subset(B_{\mu}%
^{n}(1),\widetilde{g}_{\infty,\mu})$ and $B^{n}(Q_{\infty,\nu},\delta
)\subset(B_{\nu}^{n}(1),\widetilde{g}_{\infty,\nu})$ that are isometric for
some $\delta>0.$
\end{lemma}

\begin{proof}
Choose a unit speed minimal geodesic $\beta_{i}$ joining $P_{i,\mu}=\beta
_{i}(0)$ to $P_{i,\nu}=\beta_{i}(L_{i})$. Let $Q_{i}$ be the midpoint of
$\beta_{i}$ so that $r_{i}\doteqdot d_{g_{i}}\left(  Q_{i},P_{i,\mu}\right)
=d_{g_{i}}\left(  Q_{i},P_{i,\nu}\right)  $ satisfies $\left|  2r_{i}%
-d_{\infty}(P_{\mu\infty},P_{\nu\infty})\right|  \leq1/i.$ Let $\delta=\left(
2-d_{\infty}(P_{\infty,\mu},P_{\infty,\nu})\right)  /3.$ We have by Lemma
\ref{lem overlap from geodesic} that there are isometries $\iota_{i}$ between
$\left(  B^{n}\left(  Q_{i,\mu},\delta\right)  ,\tilde{g}_{i,\mu}\right)  $
and $\left(  B^{n}\left(  Q_{i,\nu},\delta\right)  ,\tilde{g}_{i,\nu}\right)
$ which satisfy $\iota_{i}\left(  Q_{i,\mu}\right)  =Q_{i,\nu}.$ By passing to
a subsequence, we get a limit isometry $\iota_{\infty}:\left(  B^{n}%
(Q_{\infty,\mu},\delta),\widetilde{g}_{\infty,\mu}\right)  \rightarrow\left(
B^{n}(Q_{\infty,\nu},\delta),\widetilde{g}_{\infty,\nu}\right)  .$
\end{proof}

We now prove a gluing lemma. We use the notations in the proof of Lemma
\ref{lem pos curv orb}.

\begin{lemma}
\label{lem piece metric}Suppose that for all $P_{\infty}\in B_{d_{\infty
}\left(  0\right)  }\left(  O_{\infty},a\right)  \subset X_{\infty}$ there is
a neighborhood of $P_{\infty}$ isometric to $\left.  \left(  D^{2}\left(
1/2\right)  \times(-1/2,1/2),^{2}\widetilde{h}_{\infty}(t)+du^{2}\right)
\right/  G_{\infty},$ where $^{2}\widetilde{h}_{\infty}(t)$ has positive
curvature, and that there is a free isometric action of $\Delta\subset
G_{\infty}^{0}$ on $D^{2}\left(  1/2\right)  $. Then there is a metric on
$\left(  -a,a\right)  \times S^{1}$ with positive curvature such that the
diameter is at least $2a$.
\end{lemma}

\begin{proof}
For two different points $P_{\infty,\mu}$ and $P_{\infty,\nu}$ in $X_{\infty}%
$, by Lemma \ref{lem isom metr Fuk} the two metrics $^{2}\widetilde{h}%
_{\infty,\mu}(0)+du_{\mu}^{2}$ and $^{2}\widetilde{h}_{\infty,\nu}(0)+du_{\nu
}^{2}$ are isometric over some balls. We can extend the isometry to domains
which map to $B(P_{\infty,\mu},1)\cap B(P_{\infty,\nu},1)\subset X_{\infty}$
under the quotient map. Note that we have
\begin{align*}
B(P_{\infty,\mu},1)  &  \supset\left.  \left(  D^{2}\left(  1/2\right)
\times(-1/2,1/2),^{2}\widetilde{h}_{\infty,\mu}(0)+du_{\mu}^{2}\right)
\right/  G_{\infty,\mu},\\
B(P_{\infty,\nu},1)  &  \supset\left.  \left(  D^{2}\left(  1/2\right)
\times(-1/2,1/2),^{2}\widetilde{h}_{\infty,\nu}(0)+du_{\nu}^{2}\right)
\right/  G_{\infty,\nu},
\end{align*}
and $G_{\infty,\mu}^{0}=\Delta_{\mu},G_{\infty,\nu}^{0}=\Delta_{\nu}$. Since
$\Delta_{\mu}$ and $\Delta_{\nu}$ give nontrivial Killing vector fields of
$^{2}\widetilde{h}_{\infty,\mu}(0)$ and $^{2}\widetilde{h}_{\infty,\nu}(0)$
respectively, by Lemma \ref{lem coord on surf}(i) we can write
\begin{align*}
^{2}\widetilde{h}_{\infty,\mu}(0)  &  =dr_{\mu}^{2}+f_{\mu}\left(  r_{\mu
}\right)  ^{2}d\theta_{\mu}^{2}\\
^{2}\widetilde{h}_{\infty,\nu}(0)  &  =dr_{\nu}^{2}+f_{\nu}\left(  r_{\nu
}\right)  ^{2}d\theta_{\nu}^{2}%
\end{align*}
such that $f_{\mu}$ and $f_{\nu}$ are never zero and the Killing fields are
proportional to $\frac{\partial}{\partial\theta_{\mu}}$ and $\frac{\partial
}{\partial\theta_{\nu}}.$ By Lemma \ref{lem coord on surf}(ii) different
$f_{\mu}$ and $f_{\nu}$ can be pieced together to form $f$ which restricts to
$f_{\mu}$ and $f_{\nu}$.

Fix a $P_{\infty}\in X_{\infty}.$ We can define a curve in $X_{\infty}$ to be
the image of the coordinate curve $\left\{  (r,0)\in V_{P_{\infty}}\right\}
.$ Then we can extend this curve from both ends. Let $P_{1\infty}$ be an end
point, then we can use the isometry in Lemma \ref{lem isom metr Fuk} to get a
short curve which is part of $\left\{  (r,\theta_{0})\in V_{P_{1\infty}%
}\right\}  $ for some fixed $\theta_{0}$. Then we use the image of $\left\{
(r,\theta_{0})\in V_{P_{1\infty}}\right\}  $ to extend the image of curve
$\left\{  (r,0)\in V_{P_{\infty}}\right\}  $. In this way we get a curve
$\gamma\left(  s\right)  ,$ $s\in\left(  -a,a\right)  $ of length $2a$ in
$X_{\infty}$.

We can then piece together $f_{\gamma\left(  s_{i}\right)  }\left(  r\right)
$ associated to different points $\gamma(s_{i})$ to form a smooth function
$f(r),$ $r\in\left(  -a,a\right)  .$ Then $dr^{2}+f\left(  r\right)
^{2}d\theta^{2}$ defines a metric on $\left(  r,\theta\right)  \in\left(
-a,a\right)  \times S^{1}$ with positive curvature since $f_{\gamma\left(
s_{i}\right)  }\left(  r\right)  >0$ and $dr^{2}+f_{\gamma\left(
s_{i}\right)  }\left(  r\right)  ^{2}d\theta^{2}$ has positive curvature.
\end{proof}

\begin{proposition}
$X_{\infty}$ has at most one singular point of type $D^{2}/\mathbb{Z}_{p}$ or
$D^{2}/\mathbb{D}_{2p}$ where $p>1$.
\end{proposition}

More generally, we have:

\begin{proposition}
Let $(X^{n},g)$ be a complete, noncompact Riemannian orbifold with singular
points of type $D^{n}/\Gamma^{\prime}$ or $D^{n}/\left(  \Gamma^{\prime
}\rtimes\mathbb{Z}_{2}\right)  ,$ where the only fixed point of the action
$\Gamma^{\prime}$ is $\vec{0}$, $\mathbb{Z}_{2}$ acts by reflection with
respect to a hyperplane, and $\Gamma^{\prime}<\Gamma^{\prime}\rtimes
\mathbb{Z}_{2}$ is a normal subgroup. Suppose the sectional curvature
$K_{g}>0$ on $X^{n}$. Then $X^{n}$ has at most one singular point with rank
$\left\vert \Gamma^{\prime}\right\vert >1$.
\end{proposition}

\begin{proof}
We will use $d$ to denote the distance induced by $g.$ The proof is by
contradiction. Suppose the proposition is not true; let $p$ and $q$ be two
singular points, and let $\gamma$ be a unit-speed minimal geodesic from
$p=\gamma(0)$ to $q=\gamma(L)$. Without loss of generality we may assume that
$\gamma(t)$ is a smooth point for all $t\in(0,L)$ or of type $D^{n}%
/\mathbb{Z}_{2}$ for all $t\in(0,L)$, otherwise we may take two consecutive
singular points on $\gamma$ be to $p$ and $q$.

Let $\rho>0$ be a constant. There is an $\varepsilon>0$ depending on $L+\rho$
and $\rho$ so that all sectional curvatures $\kappa$ on the ball of radius
$L+\rho$ centered $p$ have $\kappa\geq\varepsilon/\rho^{2}$. Let $z$ be a
point with $d(p,z)\geq r_{0}$ where $r_{0}$ will be chosen very large (to be
specified later). The existence of $z$ follows from the completeness and
noncompactness of $X$. Let $\eta$ be a unit-speed minimal geodesic from $z$ to
$\gamma([0,L])$, say $\gamma(t_{0})=\eta(0)$ and $z$ $=\eta(r+\rho)$ for some
$r>0$.

We claim that $\gamma^{\prime}(t_{0})$ is perpendicular to $\eta^{\prime}(0)$.
Note that if $\gamma(t_{0})$ is a singular point the perpendicularity means
that there is a lift $\widetilde{\gamma}^{\prime}(t_{0})$\ of $\gamma^{\prime
}(t_{0})$ and a lift $\widetilde{\eta}^{\prime}(0)$ of $\eta^{\prime}(0)$ in
the uniformized tangent space $\widetilde{T}_{\gamma(t_{0})}M$ of the tangent
cone $T_{\gamma(t_{0})}M$ such that $\widetilde{\gamma}^{\prime}(t_{0})$ is
perpendicular to $\widetilde{\eta}^{\prime}(0)$. If $t_{0}\in(0,L)$ the claim
is clearly true since otherwise we can move $t$ to left or right of $t_{0}$ to
shorten $\eta.$ If $t_{0}=0$ (a similar argument holds for $t_{0}=L$), fix a
lifting $\widetilde{\eta}^{\prime}(0)$ of $\eta^{\prime}(0)$ and consider all
the lifting of $\gamma^{\prime}(t_{0})$ in $\widetilde{T}_{\gamma(t_{0})}M$
moved by $\Gamma^{\prime}$. Since $\gamma(t_{0})$ is singular point with rank
$\left|  \Gamma^{\prime}\right|  >1$, the sum of these lifts is zero, hence at
least one of the lifts must make an angle less than or equal to $\pi/2$. But
it cannot be less than $\pi/2$, otherwise we can shorten the distance between
$z$ and $\gamma([0,L])$. The claim is proved.

By our choice of $\varepsilon$ all sectional curvatures are greater than or
equal to $\varepsilon/\rho^{2}$ along $\eta$ for a distance $\rho$ from
$\eta(0)$. Let $\zeta_{s},$ $s\in(t_{0}-\delta,t_{0}+\delta)$ be minimal
geodesics starting from $\gamma(s)$ and going through $z$. The standard
computation shows the second variation of the length of $\zeta_{s}$ at
$s=t_{0}$ is strictly negative. Indeed let $Z=\gamma^{\prime}(t_{0})$ be the
unit tangent vector to $\gamma$ and extend $Z$ to a vector field on $\eta$ by
parallel translation. Choose a function $\varphi$ to be identically 1 within
distance $\rho$ of $\eta(0)$ along $\eta$ and then to drop linearly to zero
with slope $1/r$. The second variation of the length of $\zeta_{s}$ in the
direction $\varphi Z$ is
\[
I(\varphi Z,\varphi Z)=\int\left[  |D\varphi|^{2}-\kappa\varphi^{2}\right]
ds
\]
where $\kappa=Rm(\eta^{\prime},Z,\eta^{\prime},Z)$ is the sectional curvature.
Considering the separate contributions from the part of $\eta$ within $\rho$
of $\eta(0)$ and the part beyond we get%
\[
I(\varphi Z,\varphi Z)\leq-\varepsilon/\rho+1/r.
\]
We can take $r_{0}\leq r+\rho$ large, and hence $r$ large, enough such that
\[
1/r\leq\varepsilon/2\rho
\]
and so we get
\[
I(\varphi Z,\varphi Z)\leq-\varepsilon/2\rho.
\]
The second variation is strictly negative. The first variation of the length
of $\zeta_{s}$ at $s=t_{0}$ is zero because $\gamma^{\prime}(t_{0})$ is
perpendicular to $\eta^{\prime}(0)$. But now we see $\gamma(t_{0})$ is not the
closest point on $\gamma$ to $z$, which is a contradiction. The proposition is proved.
\end{proof}

It is well known that a two dimensional noncompact orbifold is a good
orbifold, i.e., its universal cover is a smooth manifold. We have proved the following.

\begin{theorem}
There is a universal cover $\left(  \widehat{X}_{\infty},\widehat{h}_{\infty
}(t)\right)  $ of the two dimensional limit $\left(  X_{\infty},d_{\infty
}(t)\right)  $ which is a Riemannian manifold with positive curvature for each
$t$ and is diffeomorphic to $\mathbb{R}^{2}$.
\end{theorem}

Recall that around any smooth point $P_{\infty}\in X_{\infty}$, we have a
neighborhood $\widehat{U}_{P_{\infty}}$ in $\widehat{X}_{\infty}$ such that
\[
\left(  \widehat{U}_{P_{\infty}},\widehat{h}_{\infty}(0)\right)  \cong\left.
\left(  D^{2}\left(  1/2\right)  \times(-1/2,1/2),dr^{2}+f_{P_{\infty}}\left(
r\right)  ^{2}d\theta^{2}+du^{2}\right)  \right/  \Delta,
\]
where the Killing vector field of the $\Delta$ action is $K=\left(
a\frac{\partial}{\partial\theta},b\frac{\partial}{\partial u}\right)  $ for
some constants $a$ and $b\neq0$ independent of the choice of $P_{\infty}$.
Since by construction the action of $\Delta$ preserves the metric
$^{2}\widetilde{h}_{\infty}(t)+du^{2}$ on $D^{2}\left(  1/2\right)
\times(-1/2,1/2)$ for all $t,$ there is a function $k(t)>0$ and functions
$f_{P_{\infty}}\left(  r,t\right)  $ such that $k(0)=1,$ $f_{P_{\infty}%
}\left(  r,0\right)  =f_{P_{\infty}}\left(  r\right)  $ and
\[
\left(  \widehat{U}_{P_{\infty}},\widehat{h}_{\infty}(t)\right)  \cong\left.
\left(  D^{2}\left(  1/2\right)  \times(-1/2,1/2),\text{ }k(t)^{2}%
dr^{2}+f_{P_{\infty}}\left(  r,t\right)  ^{2}d\theta^{2}+du^{2}\right)
\right/  \Delta.
\]
By the calculation in \cite{Ch}\ or \cite{CGL} the quotient metric is%
\[
\widehat{h}_{\infty}(t)=k(t)^{2}dr^{2}+\frac{f_{P_{\infty}}\left(  r,t\right)
^{2}}{1+(\frac{a}{b})^{2}f_{P_{\infty}}\left(  r,t\right)  ^{2}}d\theta^{2}.
\]
When $a\neq0$, $\widehat{h}_{\infty}(t)$ is not necessarily a solution of the RF.

However we can construct a solution of the RF on $\mathbb{R}^{2}$ by piecing
together the RF solutions $k(t)^{2}dr^{2}+f_{P_{\infty}}\left(  r,t\right)
^{2}d\theta^{2}$ when $a\neq0$. The piecing together can be done with Lemma
\ref{lem piece metric}; we stop the extension in defining the curve
$\gamma\left(  s\right)  $ when we get to $r_{0}$ and $P_{\infty}$ with
$f_{P_{\infty}}\left(  r_{0},0\right)  =0.$ This can only happen on one end
since $X_{\infty}$ is noncompact. We call such a solution the \textbf{virtual
limit} associated with the two dimensional limit $\left(  X_{\infty}%
,d_{\infty}(t)\right)  $.

\begin{theorem}
Suppose that the dimension of the limit $\left(  X_{\infty},d_{\infty
}(t)\right)  $ is two and $a\neq0$ in the Killing vector field $K=\left(
a\frac{\partial}{\partial\theta},b\frac{\partial}{\partial u}\right)  $
generated by $\Delta$ which is defined in local Fukaya theory, then there is a
rotationally symmetric surface $\left(  \overleftrightarrow{X}_{\infty
},\overleftrightarrow{h}_{\infty}(t)\right)  $ associated to $(X_{\infty
},d_{\infty}(t))$ and $\overleftrightarrow{h}_{\infty}(t)$ is a rotationally
symmetric solution of the RF with positive curvature. $\overleftrightarrow
{X}_{\infty}$ is diffeomorphic to $\mathbb{R}^{2}.$
\end{theorem}

\subsection{\strut The virtual limit associated with a 1-dimensional limit
$X_{\infty}$}

From Proposition \ref{local model} when the dimension of $X_{\infty}$ is $1$,
$X_{\infty}$ has a manifold structure, possibly with boundary. Since
$X_{\infty}$ has unbounded diameter, it is diffeomorphic to $\mathbb{R}$ or
$[0,\infty).$ When $X_{\infty}\cong\lbrack0,\infty)$ there is a unique point
$P_{\infty}\in X_{\infty}$ which fits into case (1iii), (1iv), (2ai) or
(2aii), and all other points fit into case (1i) or (1ii). We will rule out
$X_{\infty}\cong\mathbb{R}$ next.

\begin{proposition}
$X_{\infty}$ cannot be $\mathbb{R}$.
\end{proposition}

\begin{proof}
If $X_{\infty}\cong\mathbb{R},$ every point $P_{\infty}\in X_{\infty}$ fits
into case (1i) or (1ii). Then $\Delta=\mathbb{R}_{\text{loc}}^{2}$ acts freely
and locally isometrically on $\left(  D^{2}\left(  1/2\right)  \times
(-1/2,1/2),^{2}\widetilde{h}_{\infty}(0)+du^{2}\right)  $. Note that $\Delta$
cannot act trivially on the second factor $(-1/2,1/2)$, otherwise we get a
free action of $\Delta$ by isometries on $\left(  D^{2}\left(  1/2\right)
,^{2}\widetilde{h}_{\infty}(0)\right)  ,$ which by Lemma \ref{Lem Abelian act}
would give us independent Killing fields. Hence using Lemma \ref{lem coord on
surf} we may assume that $\mathbb{R}_{\text{loc}}^{2}$ acts freely
on\ $D^{2}\left(  1/2\right)  \times(-1/2,1/2)$ by $\left(  r,\theta,u\right)
\mapsto\left(  r,\theta+a\tau,u+b\tau\right)  $ and $^{2}\widetilde{h}%
_{\infty}(0)=dr^{2}+f_{P_{\infty}}\left(  r\right)  ^{2}d\theta^{2}$. There is
a subgroup of $\mathbb{R}_{\text{loc}}^{2}$ which acts freely by $\left(
r,\theta,u\right)  \mapsto\left(  r,\theta+a\tau,u\right)  $. This implies
$f_{P_{\infty}}(r)>0$.

Now we can construct a complete metric on $\left(  -\infty,\infty\right)
\times S^{1}$ with positive curvature by piecing together $k(t)^{2}%
dr^{2}+f_{P_{\infty}}\left(  r,t\right)  ^{2}d\theta^{2}$ when $a\neq0$ via
Lemma \ref{lem piece metric}. Such metric on $\left(  -\infty,\infty\right)
\times S^{1}$ cannot exist as argued before. We obtain a contradiction and the
proposition is proved.
\end{proof}

When $X_{\infty}\cong\lbrack0,\infty)$ we construct a 2-dimensional virtual
limit. A simple construction follows. A more geometric construction is given
in next section.

\begin{theorem}
\label{Thm 1D virt mod}There is a rotationally symmetric surface $\left(
\overleftrightarrow{X}_{\infty},\overleftrightarrow{h}_{\infty}(t)\right)  $
diffeomorphic to $\mathbb{R}^{2}$ associated to $X_{\infty}\cong\left(
\lbrack0,\infty),h_{\infty}(t)\right)  $ and $\overleftrightarrow{h}_{\infty
}(t)$ is a solution of the RF with positive curvature.
\end{theorem}

\begin{proof}
For every point $P_{\infty}\in(0,\infty)$, we have $^{2}\widetilde{h}_{\infty
}(t)=k(t)^{2}dr^{2}+f_{P_{\infty}}\left(  r,t\right)  ^{2}d\theta^{2}.$ We can
construct a metric $k(t)^{2}dr^{2}+f\left(  r,t\right)  ^{2}d\theta^{2}$ with
positive curvature on $\left(  0,\infty\right)  \times S^{1}$ such that the
projection $\left(  0,\infty\right)  \times S^{1}\rightarrow(0,\infty)\cong
X_{\infty}-\left\{  O_{\infty}\right\}  $ is a submersion. This can be done by
piecing together using Lemma \ref{lem piece metric}. To see what happens at
$0\in X_{\infty}$, we choose $P_{\infty}=0$ and assume the corresponding
coordinate $r(P_{\infty})=0$. Then $0$ fits into case (1iii), (1iv), (2ai), or (2aii).

If $0$ fits into either case (1iii) or (1iv), then $f_{P_{\infty}}\left(
0,0\right)  \neq0$ from the fact that $\Delta=\mathbb{R}_{\text{loc}}^{2}$
acts freely and $f_{P_{\infty}}\left(  -r,0\right)  =f_{P_{\infty}}\left(
r,0\right)  $ from the $\mathbb{Z}_{2}$\ reflection. We can extend the metric
$k(0)^{2}dr^{2}+f\left(  r,0\right)  ^{2}d\theta^{2}$ using the $\mathbb{Z}%
_{2}$\ reflection to get a\ complete metric on $\left(  -\infty,\infty\right)
\times S^{1}$ with positive curvature. Such metric on $\left(  -\infty
,\infty\right)  \times S^{1}$ cannot exist as argued before. Hence case (1iii)
or (1iv) can not happen.

If $0$ fits into either case (2ai), (2aii), then $f_{P_{\infty}}\left(
0,0\right)  =0$ from the fact that the action of $\Gamma$ has fixed point $0$.
Hence $f\left(  0,t\right)  =0$ and $k(t)^{2}dr^{2}+f\left(  r,t\right)
^{2}d\theta^{2}$ extend to a smooth metric on $\left(  [0,\infty)\times
S^{1}\right)  /\left(  \left\{  0\right\}  \times S^{1}\right)  \cong
\mathbb{R}^{2}$ which is our $\overleftrightarrow{X}_{\infty}$. It is clear
that $k(t)^{2}dr^{2}+f\left(  r,t\right)  ^{2}d\theta^{2}$ is a rotationally
symmetric solution of the RF with positive curvature.
\end{proof}

\subsection{A geometric construction}

We now give a geometric, but more technical description of how to construct
the virtual limit of the one-dimensional $X_{\infty}$.

\textbf{Step 1}: Create framework. Consider the covering $\left\{
U_{\infty,k}\right\}  _{k\in\mathbb{N}\cup\left\{  0\right\}  }$ of
$X_{\infty}$ defined by $U_{\infty,0}=[0,3)$ and $U_{\infty,k}=(4k-3,4k+3)$
for $k\in\mathbb{N}$. It is clear that $U_{\infty,k}\cap U_{\infty,\ell}%
\neq\varnothing$ if and only if $|k-\ell|\leq1.$ Let $P_{\infty,k}=4k,$ which
are the centers of $U_{k}.$ A word about notation: $B_{i,k}^{3}\left(
3\right)  \subset\mathbb{R}^{3}$ all denote the same 3-ball of radius $3,$
where the indices $i\in\mathbb{N}\cup\left\{  \infty\right\}  $ and
$k\in\mathbb{N}\cup\left\{  0\right\}  $ are just to remind the reader of the
dependence on the term in the sequence and the point in $X_{\infty}.$ We also
have frames at $P_{i,k},$ orthonormal with respect to $g_{i}\left(  0\right)
,$ which we use to identify the tangent spaces $T_{P_{i,k}}M_{i}$ with
$\mathbb{R}^{3}.$

\textbf{Step 2}: Glue open sets in the sequence. By the definition of GH
convergence, passing to a subsequence there exist maps $\varphi_{i}\left(
t\right)  :\left(  M_{i},d_{g_{i}\left(  t\right)  },O_{i}\right)
\rightarrow\left(  X_{\infty},d_{\infty}\left(  t\right)  ,O_{\infty}\right)
$ and $\psi_{i}\left(  t\right)  :\left(  X_{\infty},d_{\infty}\left(
t\right)  ,O_{\infty}\right)  \rightarrow\left(  M_{i},d_{g_{i}\left(
t\right)  },O_{i}\right)  $ for each $t\in\left[  \beta,\phi\right]  $ which
are $1/i$-pointed GH approximations. By choosing a countable dense set of
times in $\left[  \beta,\psi\right]  $ and passing to a subsequence we can
find time-independent maps $\varphi_{i}:\left(  M_{i},d_{g_{i}\left(
t\right)  },O_{i}\right)  \rightarrow\left(  X_{\infty},d_{\infty}\left(
t\right)  ,O_{\infty}\right)  $ and $\psi_{i}:\left(  X_{\infty},d_{\infty
}\left(  t\right)  ,O_{\infty}\right)  \rightarrow\left(  M_{i},d_{g_{i}%
\left(  t\right)  },O_{i}\right)  $ which are $1/i$-pointed GH approximations
for all $t\in\left[  \beta,\psi\right]  .$ Let $P_{i,k}\doteqdot\psi
_{i}\left(  P_{\infty,k}\right)  $ for $k<i/4.$ Choose a minimal geodesic
$\beta_{i,k}$ with respect to $g_{i}\left(  0\right)  $ joining $P_{i,k}$ to
$P_{i,k+1}$ for $k+1<i/4.$ Let $Q_{i,k}$ be the midpoint of $\beta_{i,k}$ so
that $r_{i,k}\doteqdot d_{g_{i}\left(  0\right)  }\left(  Q_{i,k}%
,P_{i,k}\right)  =d_{g_{i}\left(  0\right)  }\left(  Q_{i,k},P_{i,k+1}\right)
$ satisfies $\left|  2r_{i,k}-4\right|  \leq1/i.$ For $i$ large enough we have
$B_{g_{i}\left(  0\right)  }\left(  Q_{i,k},1/2\right)  \subset B_{g_{i}%
\left(  0\right)  }\left(  P_{i,k},3\right)  \cap B_{g_{i}\left(  0\right)
}\left(  P_{i,k+1},3\right)  .$ For each $i\in\mathbb{N}$ and $t\in\left[
\beta,\psi\right]  ,$ let $\tilde{g}_{i,k}\left(  t\right)  =\exp_{P_{i,k}%
}^{\ast}g_{i}\left(  t\right)  $ (recall the definition from \S\ref{sec 3.2
limit metric}) and consider the 3-balls of radius three $\left(
B_{i,k}\left(  3\right)  ,\tilde{g}_{i,k}\left(  t\right)  \right)  $
isometrically covering $\left(  B_{g_{i}\left(  0\right)  }\left(
P_{i,k},3\right)  ,g_{i}\left(  t\right)  \right)  $ by the exponential map
$\exp_{P_{i,k}}$ for each $k\in\mathbb{N}\cup\left\{  0\right\}  .$

We can use $\beta_{i,k}$ to identify $\left(  B_{i,k}\left(  3\right)
,\tilde{g}_{i,k}\left(  0\right)  \right)  $ and $\left(  B_{i,k+1}\left(
3\right)  ,\tilde{g}_{i,k+1}\left(  0\right)  \right)  $ by an isometry
$\iota_{i,k}$ on the overlap regions (which we define below). Moreover,
$\iota_{i,k}$ is an isometry with respect to the metrics $\tilde{g}%
_{i,k}\left(  t\right)  $ and $\tilde{g}_{i,k+1}\left(  t\right)  $ for all
$t\in\left[  \beta,\psi\right]  .$

As in Lemma \ref{lem overlap from geodesic} there are balls around points
$\tilde{Q}_{i,k}$ and $\overset{\approx}{Q}_{i,k}$ respectively in
$B_{i,k}\left(  3\right)  $ and $B_{i,k+1}\left(  3\right)  $ which are
isometric and $\exp_{P_{i,k}}\tilde{Q}_{i,k}=\exp_{P_{i,k+1}}\overset{\approx
}{Q}_{i,k}=Q_{i,k}$. Let $\tilde{U}_{i,k}$ be the connected component of
\[
\left(  \exp_{P_{i,k}}\right)  ^{-1}\left(  B_{g_{i}\left(  0\right)  }\left(
P_{i,k},3\right)  \cap B_{g_{i}\left(  0\right)  }\left(  P_{i,k+1},3\right)
\right)
\]
containing $\tilde{Q}_{i,k}$ and $\overset{\approx}{U}_{i,k}$ be the connected
component of
\[
\left(  \exp_{P_{i,k+1}}\right)  ^{-1}\left(  B_{g_{i}\left(  0\right)
}\left(  P_{i,k},3\right)  \cap B_{g_{i}\left(  0\right)  }\left(
P_{i,k+1},3\right)  \right)
\]
containing $\overset{\approx}{Q}_{i,k}.$ By local covering space theory, there
exists a unique diffeomorphism $\iota_{i,k}:\left(  \tilde{U}_{i,k},\tilde
{Q}_{i,k}\right)  \rightarrow\left(  \overset{\approx}{U}_{i,k},\overset
{\approx}{Q}_{i,k}\right)  $ such that $\left.  \exp_{P_{i,k}}\right|
_{\tilde{U}_{i,k}}=\exp_{P_{i,k+1}}\circ\,\iota_{i,k}$ which extends the
isometry from Lemma \ref{lem overlap from geodesic} mentioned above. The map
$\iota_{i,k}$ is an isometry from $\tilde{g}_{i,k}\left(  t\right)  $ to
$\tilde{g}_{i,k+1}\left(  t\right)  $ for all $t\in\left[  \beta,\psi\right]
.$ By passing to a subsequence, we get a limit isometry $\iota_{\infty,k}$ of
the overlap regions $\tilde{U}_{\infty,k}$ of $\left(  B_{\infty,k}\left(
3\right)  ,\tilde{g}_{\infty,k}\left(  t\right)  \right)  $ and $\overset
{\approx}{U}_{\infty,k}$ of $\left(  B_{\infty,k+1}\left(  3\right)
,\tilde{g}_{\infty,k+1}\left(  t\right)  \right)  $ for all $k\in
\mathbb{N}\cup\left\{  0\right\}  $ and $t\in\left[  \beta,\psi\right]  .$

\textbf{Step 3}: Find good coordinates. From $B_{\infty,0}\left(  3\right)  $
we get rotationally symmetric metrics $h_{\infty,0}\left(  t\right)  $ on a
2-disk $\Sigma_{\infty,0}^{2}=D^{2}\left(  r_{\infty,0}\right)  $ for some
$r_{\infty,0}\geq2$. We obtain these metrics by considering the surface slice
passing through $\vec{0}$ (recall that metrically the ball is locally a
surface product with a line). Alternatively, these metrics can be obtained
from the isometric embedding of $D^{2}\left(  r_{\infty,0}\right)
\times\left(  -\varepsilon,\varepsilon\right)  $ into $\left(  B_{\infty
,0}\left(  3\right)  ,\tilde{g}_{\infty,k}\left(  t\right)  \right)  $ by
quotienting out the interval direction. For $t=0,$ letting $\varepsilon
\rightarrow0$ allows us to take $r_{\infty,0}\rightarrow3.$ For $k\in
\mathbb{N}$ from $B_{\infty,k}\left(  3\right)  $ we get metrics $\left(
\Sigma_{\infty,k}^{2},h_{\infty,k}\left(  t\right)  \right)  $ which are
locally warped products. Let $O_{\infty,k}$ denote the point in $\Sigma
_{\infty,k}^{2}$ corresponding to $\vec{0}\in B_{\infty,k}\left(  3\right)  .$
Then $\overline{B\left(  O_{\infty,k},r\right)  }\subset\Sigma_{\infty,k}^{2}$
is compact for all $r<3.$ The warped product metric $h_{\infty,1}\left(
0\right)  $ on $\Sigma_{\infty,k}^{2}$ induces oriented local coordinates
$\left(  r,\theta\right)  $ such that $h_{\infty,k}\left(  0\right)
=dr^{2}+f_{k}\left(  r\right)  ^{2}d\theta^{2}.$ The coordinate $r$ is
uniquely determined up an additive constant and $\theta$ is uniquely
determined by $\theta\left(  O_{\infty,k}\right)  =0$ by Lemma \ref{lem coord
on surf}. For every $\varepsilon>0,$ there exists $\delta>0$ such that%
\[
\left(  r,\theta\right)  :\left(  -3+\varepsilon,3-\varepsilon\right)
\times\left(  -\delta,\delta\right)  \rightarrow\Sigma_{\infty,k}^{2}%
\]
is an embedding.

\textbf{Step 4}: Form metrics on disks. This step corresponds to Lemma
\ref{lem piece metric}. We have constructed an overlap map from $\Sigma
_{\infty,0}^{2}$ to $\Sigma_{\infty,1}^{2}.$ Then we extend $S_{\infty,1}%
^{2}=\Sigma_{\infty,0}^{2}\cup\Sigma_{\infty,1}^{2}/\sim$ (alternately we
could have defined $S_{\infty,1}^{2}=\Sigma_{\infty,0}^{2}\cup\left[  \left(
-3+\varepsilon,3-\varepsilon\right)  \times\left(  -\delta,\delta\right)
\right]  /\sim$), where the identification $\sim$ is defined by the isometry
$\iota_{\infty,0},$ to rotationally symmetric metrics $j_{\infty,1}\left(
t\right)  $ on a disk $\Delta_{\infty,1}^{2}=D^{2}\left(  r_{\infty,1}\right)
.$ One can see this explicitly in terms of the coordinates $\left(
r_{0},\theta_{0}\right)  $ and $\left(  r_{1},\theta_{1}\right)  $ on
$\Sigma_{\infty,0}^{2}-\left\{  O_{\infty,0}\right\}  $ and $\Sigma_{\infty
,1}^{2}$ with $\lim_{x\rightarrow O_{\infty,0}}r_{0}\left(  x\right)  =0$ and
$r_{1}\left(  O_{\infty,1}\right)  =0.$ In these coordinates, the isometry
$\iota_{\infty,0}$ is given by $r_{1}=r_{0}-4$ and $\theta_{1}=c_{1}\theta
_{0}+c_{2},$ for some $c_{1},c_{2}\in\mathbb{R},$ and maps $\left(  \left(
1+\varepsilon,3-\varepsilon\right)  \times\left(  -\delta,\delta\right)
,h_{\infty,0}\left(  t\right)  \right)  $ to $\left(  \left(  -3+\varepsilon
,-1-\varepsilon\right)  \times\left(  -\delta,\delta\right)  ,h_{\infty
,1}\left(  t\right)  \right)  .$ Under this isometry at $t=0$ we have
$f_{1}\left(  r-4\right)  =f_{0}\left(  r\right)  $ for $r\in\left(
1+\varepsilon,3-\varepsilon\right)  .$

Define the smooth function $f$ as
\[
f\left(  r\right)  =\left\{
\begin{tabular}
[c]{ll}%
$f_{0}\left(  r\right)  $ & $0\leq r\leq2$\\
$f_{1}\left(  r-4\right)  $ & $2<r\leq6$%
\end{tabular}
\right.  .
\]
This gives us a metric $dr^{2}+f\left(  r\right)  ^{2}d\theta^{2}$ on
$D^{2}\left(  6\right)  .$ We can also construct this space by starting on
$S_{\infty,1}^{2}$ and then using the symmetry to extend to $D^{2}\left(
6\right)  $ via
\[
D^{2}\left(  6\right)  =S_{\infty,1}^{2}\times S^{1}/\sim,
\]
where $\left(  r,\theta,\phi\right)  \sim\left(  r,\theta+a,\phi+a\right)  $
for all $a\in S^{1}$ and the metric is extended so the translation of the
circle (denoted by $a$ above) is isometric. We can continue this process to
$D^{2}\left(  k+2\right)  $ for all $k.$ This produces a complete metric
$dr^{2}+f\left(  r\right)  ^{2}d\theta^{2}$ on $\mathbb{R}^{2}.$ This
construction works for any $t.$ We take as the virtual limit the warped
product metric $\overleftrightarrow{h}(t)$ with $\overleftrightarrow
{h}(0)=dr^{2}+f\left(  r\right)  ^{2}d\theta^{2}$ on $\mathbb{R}^{2}.$ Note
that every point has a neighborhood isometric to $\left(  \Sigma_{\infty
,k}^{2},h_{\infty,k}\left(  t\right)  \right)  $ for some $k\geq0.$ Note that
after a translation of $r_{k}$ and $\theta_{k}$ and scaling of $\theta_{k}$
the coordinates $\left(  r_{k},\theta_{k}\right)  $ define isometric
immersions of $\left(  \Sigma_{\infty,k}^{2},h_{\infty,k}\left(  t\right)
\right)  $ into $\left(  \mathbb{R}^{2},\overleftrightarrow{h}(t)\right)  .$

\section{Hamilton's singularity theory in the collapsed Type IIb case}

Let $\left(  M^{3},g\left(  t\right)  \right)  ,$ $t\in\lbrack0,\infty),$ be a
complete solution of the Ricci flow with bounded curvature forming a Type IIb
singularity at time $\infty.$ That is, $\sup_{M\times\lbrack0,\infty)}t\left|
Rm\left(  x,t\right)  \right|  =\infty.$ We follow \S16 of \cite{Hamilton
Formation} in choosing a sequence of points and times to dilate about. Let
$T_{i}\rightarrow\infty$ be any sequence of times. Given any sequence
$\varepsilon_{i}\rightarrow0,$ choose $\left(  x_{i},t_{i}\right)  \in
M\times\left[  0,T_{i}\right]  $ so that $t_{i}\rightarrow\infty$ and
\[
\frac{t_{i}\left(  T_{i}-t_{i}\right)  \left|  Rm\left(  x_{i},t_{i}\right)
\right|  }{\sup_{M\times\left[  0,T_{i}\right]  }\left\{  t\left(
T_{i}-t\right)  \left|  Rm\left(  x,t\right)  \right|  \right\}  }%
\geq1-\varepsilon_{i}.
\]
This can be done because
\[
\sup_{M\times\left[  0,T_{i}\right]  }\left\{  t\left(  T_{i}-t\right)
\left|  Rm\left(  x,t\right)  \right|  \right\}  \geq\frac{T_{i}}{2}%
\sup_{M\times\left[  0,T_{i}/2\right]  }\left\{  t\left|  Rm\left(
x,t\right)  \right|  \right\}  ,
\]
and the right side goes to infinity. Let $K_{i}=\left|  Rm\left(  x_{i}%
,t_{i}\right)  \right|  .$ The solution $g_{i}\left(  t\right)  \doteqdot
K_{i}g\left(  t_{i}+t/K_{i}\right)  $ exists on the time interval
$[-\alpha_{i},\infty),$ where $\alpha_{i}=t_{i}K_{i}\rightarrow\infty$ as
$i\rightarrow\infty.$ Let $\omega_{i}=\left(  T_{i}-t_{i}\right)  K_{i}.$ Then
$\omega_{i}\rightarrow\infty.$ We have%
\[
\left|  Rm\left[  g_{i}\right]  \left(  x,t\right)  \right|  \leq\frac
{1}{1-\varepsilon_{i}}\frac{\alpha_{i}}{\alpha_{i}+t}\frac{\omega_{i}}%
{\omega_{i}-t}%
\]
for all $x\in M$ and $t\in\left[  -\alpha_{i},\omega_{i}\right]  ,$ as well as
$\left|  Rm\left[  g_{i}\right]  \left(  x_{i},0\right)  \right|  =1.$ Note
that $g_{i}\left(  t\right)  $ have ANSC by \S\ref{sec 2.3 dilate}.

Now consider the sequence $\left(  M^{3},g_{i}\left(  t\right)  ,x_{i}\right)
,$ $t\in\lbrack-\alpha_{i},\omega_{i}],$ and assume that the diameters tend to
infinity and the sequence collapses. Recall that this implies the origins are
split-like. The virtual limit $\left(  \overleftrightarrow{X}_{\infty
},\overleftrightarrow{h}_{\infty}(t)\right)  $ constructed is a complete
solution on a surface diffeomorphic to $\mathbb{R}^{2}$ with positive, bounded
curvature and exists for $t\in\left(  -\infty,\infty\right)  .$ The base point
$\overleftrightarrow{x}_{\infty}\in\overleftrightarrow{X}_{\infty}$ satisfies
\[
\left|  Rm\left[  \overleftrightarrow{h}_{\infty}(t)\right]  \left(
\overleftrightarrow{x}_{\infty},0\right)  \right|  =\sup_{\overleftrightarrow
{X}_{\infty}\times\left(  -\infty,\infty\right)  }\left|  Rm\left[
\overleftrightarrow{h}_{\infty}(t)\right]  \left(  x,t\right)  \right|  .
\]
The reason is that any other point $\left(  \overleftrightarrow{y}_{\infty
},0\right)  \in\overleftrightarrow{X}_{\infty}\times\left(  -\infty
,\infty\right)  $ corresponds to a sequence $\left(  y_{i},0\right)  \in
M\times\left[  \alpha_{i},\omega_{i}\right]  $ endowed with the metric
$g_{i}\left(  0\right)  .$ The curvatures at $\left(  y_{i},0\right)  $ are
almost smaller than the curvatures at $\left(  x_{i},0\right)  ,$ and hence
the curvatures of the limit of the covering geometries of the $y_{i}$ are less
than those of the $x_{i}.$ The curvature of the virtual limit is the same as
the curvature in the surface direction of the limit of the covering geometries
(which splits as a product), and hence the curvature at $\left(
\overleftrightarrow{y}_{\infty},0\right)  $ is less than the curvature at
$\left(  \overleftrightarrow{x}_{\infty},0\right)  .$

Hence, by Hamilton's result that eternal solutions are steady solitons
\cite{Hamilton Eternal}, which uses his matrix Harnack estimate, $\left(
\overleftrightarrow{X}_{\infty},\overleftrightarrow{h}_{\infty}(t)\right)  $
is isometric to a cigar soliton solution. Here we have used the fact that
since $\overleftrightarrow{X}_{\infty}$ is 2-dimensional, the maximum of
$\left|  Rm\right|  $ being attained at $\left(  \overleftrightarrow
{x}_{\infty},0\right)  $ is the same as the maximum of $R$ being attained at
$\left(  \overleftrightarrow{x}_{\infty},0\right)  .$ This last condition is
what is assumed in Hamilton's eternal solutions result.

\end{document}